%
% August, 2009,
% Paper Id: math.DG/0404326
% Paper pw: jbxwj
%

\input amstex
\documentstyle{amsppt}

\magnification=1200
 \vsize19.5cm
 \hsize13.8cm
\TagsOnLeft
 \pageno=1
%\nopagenumbers

%\def\endbox{\hbox{\vrule height8pt width5pt depth2pt}}

\def\p{\partial}
\def\D{\nabla}
\def\noo{\noindent}
\def\eps{\varepsilon}
\def\lam{\lambda}
\def\R{\bold R}
\def\th{\theta}
\def\wtt{\widetilde}
\def\back{\backslash}
\def\det{\text{det}}

\def\lan{\langle}
\def\ran{\rangle}

\def\arctg{\text{arctg}}
\def\dist{\text{dist}}
\def\phi{\varphi}
\def\ol{\overline}
\def\ul{\underline}
\def\Om{\Omega}

\def\bom{\ol \Om}
\def\pom{\p\Om}
\def\Ga{\Gamma}

\def\div{\text{div}}
\def\hx{\hat x}
\def\wx{\wtt x}
\def\M{\Cal M}
\def\L{\Cal L}

\def\F{\Cal F}

\NoRunningHeads
\nologo

\voffset=-5pt

\topmatter

\title{Convex solutions\\ to the mean curvature flow
}\endtitle

\author{Xu-Jia Wang}\endauthor

\affil{Centre for Mathematics and Its Applications\\
       The Australian National University\\
      Canberra, ACT 0200\\
      Australia\\
 }\endaffil

\address{Centre for Mathematics and Its Applications,
       Australian National University, \newline
       Canberra, ACT 0200,
       Australia}\endaddress

\email{wang\@maths.anu.edu.au}\endemail

\thanks{This work was supported
        by the Australian Research Council.}\endthanks

%\keywords {mean curvature flow.}\endkeywords


\abstract {In this paper we study the classification of ancient
convex solutions to the mean curvature flow in $\R^{n+1}$. An open
problem related to the classification of type I\!I singularities
is whether a convex translating solution is $k$-rotationally
symmetric for some integer $2\le k\le n$, namely whether its level
set is a sphere or cylinder $S^{k-1}\times \R^{n-k}$.  In this
paper we give an affirmative answer for entire solutions in
dimension 2. In high dimensions we prove that there exist
non-rotationally symmetric, entire convex translating solutions,
but the blow-down in space of any entire convex translating
solution is $k$-rotationally symmetric. We also prove that the
blow-down in space-time of an ancient convex solution which sweeps
the whole space $\R^{n+1}$ is a shrinking sphere or cylinder.
 }\endabstract

\endtopmatter

\vskip-15pt

\document

\baselineskip=14.0pt
\parskip=2.5pt

\centerline{\bf 1. Introduction}

\vskip10pt

Convex solutions arise in the study of singularities of the mean
curvature flow. To study the geometric behavior at singularities one
needs to classify such solutions. In this paper\footnote{
       This is a revised version of the paper arXiv:math.DG/0404326.
       All results and the ideas of proofs are the same. The main change
       is the proof of the growth estimate (1.5) in Section 2.
       In this new version we divide it into two parts.
       We first prove it for dimension 2, then prove it for high dimensions.} we study the
classification, or more precisely the geometric asymptotic behavior,
of general ancient convex solutions, including the convex
translating solutions arising at type I\!I singularities [14,15] and
the ancient convex solutions arising at general singularities [28].

It was proved by Huisken-Sinestrari [14, 15] that if $\M$ is a
mean convex flow, namely a mean curvature flow with mean convex
solution, in the Euclidean space $\R^{n+1}$, then the limit flow
obtained by a proper blow-up procedure near type I\!I singular
points is a convex translating solution (also called soliton),
that is in an appropriate coordinate system a mean curvature flow
of the form $\M'=\{(x, u(x)+t)\in \R^{n+1}:\ x\in\R^n, t\in\R\}$,
where $u$ is a complete convex solution to the mean curvature
equation
$$\div(\frac{Du}{\sqrt{1+|Du|^2}})
           =\frac {1}{\sqrt{1+|Du|^2}} . \tag 1.1 $$

Translating solutions play a similar role for the investigation of
asymptotic behavior of type I\!I singularities as self-similar
solutions for type I singularities. It is known that a convex
self-similar solution must be a shrinking sphere or cylinder [13].
For convex translating solutions there is a well known conjecture
among researches in this area, which is explicitly formulated, for
example in [28], which asserts that if $u$ is a complete convex
solution of (1.1), then the level sets $\{u=\text{const}\}$ are
spheres or cylinders. This Bernstein type problem attracted
attention in recent years, as it is crucial for a classification
of type I\!I singularities of the mean convex flow. In this paper
we prove the conjecture is true for entire solutions in dimension
two (Theorem 1.1) and false in higher dimensions (Theorem 1.2).

In this paper we also study the classification of general ancient
convex solutions to the mean curvature flow. In [28] White proved
that any limit flow to the mean convex flow in $\R^{n+1}$ for $n<7$,
or any special limit flow, namely blowup solution before first time
singularity for $n\ge 7$, is an ancient convex solution, namely at
any time the solution is a convex hypersurface. We prove that an
ancient convex solution is convex in space-time (Proposition 4.1);
and that the parabolic blow-down in space-time of any entire,
ancient convex solution, and the blow-down in space of any entire
convex translating solution, is a shrinking sphere or cylinder
(Theorem 1.3). This result corresponds to Perelman's classification
of ancient $\kappa$-noncollapsing solutions with nonnegative
sectional curvature to the 3-dimensional Ricci flow [22], see \S6.

To study the above two problems, we will consider the following
more general equation
$$\Cal L_\sigma [u] =: \sum_{i,j=1}^n
   (\delta_{ij}-\frac {u_iu_j}{\sigma +|Du|^2})u_{ij} =1,\tag 1.2$$
where $\sigma\in [0, 1]$ is a constant. If $u$ is a convex
solution of (1.2), then $u+t$, as a function of $(x, t)\in
\R^n\times\R$, is a translating solution to the flow
$$u_t=\sqrt{\sigma+|Du|^2}
            \div(\frac {Du}{\sqrt{\sigma+|Du|^2}}).\tag 1.3$$
When $\sigma=1$, equation (1.2) is exactly the mean curvature
equation (1.1), and (1.3) is the non-parametrized mean curvature
flow. When $\sigma=0$, (1.3) is the level set flow. That is if $u$
is a solution of (1.2) with $\sigma=0$,  the level set $\{u=-t\}$,
where $-\infty<t<-\inf u$, evolves by mean curvature.

Conversely, if a family of convex hypersurfaces $\M=\{\M_t\}$,
with time slice $\M_t$, evolves by mean curvature, then $\M$ can
be represented as a graph of $u$ in the space-time
$\R^{n+1}\times\R^1$ with $x_{n+2}=-t$, and the function $u$
satisfies (1.2) with $\sigma=0$. We will show that the function
$u$ itself is convex (Proposition 4.1). Therefore for both
problems it suffices to study the classification of convex
solutions to equation (1.2).

We say a solution to the mean curvature flow is {\it ancient} if
it exists from time $-\infty$. We say a solution $u$ of (1.2) is
{\it complete} if its graph is a complete hypersurface in
$\R^{n+1}$, and $u$ is an {\it entire} solution if it is defined
in the whole space $\R^n$. Accordingly an ancient convex solution
$\M$ to the mean curvature flow in $\R^{n+1}$ is an entire
solution if $\M$ is an entire graph in space-time
$\R^{n+1}\times\R^1$, which is equivalent to saying that the flow
$\M$ sweeps the whole space $\R^{n+1}$. We say $u$ is {\it
$k$-rotationally symmetric} if there exists an integer $1\le k\le
n$ such that in an appropriate coordinate system, $u$ is
rotationally symmetric with respect to $x_1, \cdots, x_k$ and is
independent of $x_{k+1},\cdots, x_n$. Therefore a function $u$ is
$k$-rotationally symmetric if and only if its level sets are
spheres ($k=n$) or cylinders ($k<n$). For other related
terminologies we refer the reader to [15, 28]. For any $1\le k\le
n$, there is a $k$-rotationally symmetric convex solution of
(1.1), which is unique up to orthogonal transformations. When
$n=1$, the unique complete convex solution of (1.1) is the
$``$grim reaper$"$, given by $u(x)=\log\sec x_1$. To exclude
hyperplanes in this paper we always consider convex solution with
positive mean curvature.

The results in this paper can be summarized in the following
theorems.

\proclaim{Theorem 1.1} If $n=2$, then any entire convex solution to
(1.2) must be rotationally symmetric in an appropriate coordinate
system.
\endproclaim

From Theorem 1.1 we obtain

\proclaim{Corollary 1.1} A convex translating solution to the mean
curvature flow must be rotationally symmetric if it is a limit
flow to a mean convex flow in $\R^3$.
\endproclaim

\proclaim{Theorem 1.2} For any dimension $n\ge 2$ and $1\le k\le
n$, there exist complete convex solutions, defined in strip
regions, to equation (1.2) which are not $k$-rotationally
symmetric. If $n\ge 3$, there exist entire convex solutions to
(1.2) which are not $k$-rotationally symmetric.
\endproclaim

Theorems 1.1 and 1.2 reflect a typical phenomenon, namely the
Bernstein theorem is in general true in low dimensions and false
in higher dimensions. See [26] for a brief discussion.

\proclaim{Theorem 1.3} Let $u$ be an entire convex solution of
(1.2). Let $u_h(x)=h^{-1}u\sqrt h\, x)$. Then there is an integer
$2\le k\le n$ such that after a rotation of the coordinate system
for each $h$, $u_h$ converges to
$$\eta_k (x) =\frac 1{2(k-1)}\sum_{i=1}^k x_i^2.\tag 1.4 $$
\endproclaim

The case $\sigma=0$ of Theorems 1.1- 1.3 describes the geometry of
ancient convex solutions to the mean curvature flow, while the case
$\sigma=1$ of Theorems 1.1-1.3 resolves the problem on convex
translating solutions.  Note that if $u$ is a convex solution which
is not defined in the whole space, then $u$ is defined in a convex
strip region (Corollaries 2.1 and 2.2), and it cannot be a blowup
solution to the mean convex flow in general (Corollary 6.1). We also
remark that in Theorem 1.3 we did not rule out the possibility that
the axis of the cylinder-like level set $\{u_h=1\}$ may rotate
slowly as $h\to\infty$, as the convergence $u_h\to \eta_k$ is
uniform only on any compact sets.

As a limit flow at the first time singularity is convex, by
Brakke's regularity theorem, a blowup sequence converges smoothly
on any compact sets to an ancient convex solution [28]. Therefore
by the above classifications one may infer that if $\M=\{\M_t\}$
is a mean convex flow in $\R^{n+1}$, $n\ge 2$, then $\M_t$
satisfies a canonical neighborhood condition, similar to the
assertion in [23] for the Ricci flow, at any point $x_t\in\M_t$
with large mean curvature before the first time singularity. In
particular if the mean curvature at $x_t$ converges to infinity,
then a proper scaling of $\M$ at $x_t$ converges along
subsequences to shrinking spheres or cylinders. See discussions in
\S6.

\vskip5pt

Our proofs of the above theorems rely heavily on the convexity of
solutions. To prove these theorems it suffices to consider the
cases $\sigma=0$ and $\sigma =1$, as for any $\sigma>0$, one can
make the transform $\hat u(x)=\frac 1\sigma u(\sqrt\sigma x)$ to
change equation (1.2) to the case $\sigma=1$. A key estimate for
the proof of Theorem 1.3 is that for any entire convex solution
$u$ of (1.2), there exists a positive constant $C$ such that
$$u(x)\le C(1+|x|^2)  \ \ \ \forall\ x\in\R^n . \tag 1.5$$
The constant $C$ depends only on $n$ and the upper bound of $u(0)$
and $|Du(0)|$. Note that (1.5) implies the compactness of the set of
entire convex solutions to (1.2), see Corollary 2.3.

By Theorem 1.3 and estimate (1.5) we have, if $n=2$,
$$C_1 |x|^2\le u(x)\le C_2|x|^2$$
for large $|x|$. Hence the case $\sigma=0$ of Theorem 1.1 follows
immediately from the asymptotic estimates in [8]. For the case
$\sigma=1$ we will prove furthermore, by an iteration argument, that
$$|u(x)-u_0(x)|=o(|x|)\ \ \ \text{as}\ \ |x|\to\infty,\tag 1.6$$
where $u_0$ is the radial solution of (1.1). We then conclude
$u=u_0$ by a Liouville type theorem of Bernstein [2], which
asserts that an entire solution $w$ to an elliptic equation in
$\R^2$ must be a constant if $|w(x)|=o(|x|)$ at infinity [24].

The proof of Theorem 1.2 is different for the cases $\sigma=0$ and
$\sigma=1$. For the case $\sigma=0$, we consider the Dirichlet
problem
$$\cases
\Cal L_0[u]  =1
        \ \ &\text{in}\ \ \Om,\\
u  = 0\ \ \ \ &\text{on}\ \ \pom, \\
\endcases \tag 1.7 $$
where $\Om$ is a bounded convex domain in $\R^n$ ($n\ge 2$). The
existence and uniqueness of viscosity solutions to (1.7) can be
found in [3,7]. We prove that there exists a sequence of bounded
convex domains $\{\Om_k\}$ such that $u_k+|\inf_{\Om_k} u_k|$,
where $u_k$ is the solution of (1.7) with $\Om=\Om_k$, converges
to a complete convex solution $u$ of $\Cal L_0[u]=1$ of which the
level set $\{x\in\R^n:\ u (x)=1\}$ is not a sphere. To prove that
$u$ is a complete convex solution we need the concavity of the
function $\log (-u)$ (Lemma 4.1).

The concavity of $\log (-u)$ is still an open problem for the mean
curvature equation (1.1). To construct a similar sequence of
solutions $(u_k)$ for equation (1.1), we use the Legendre
transform to convert the mean curvature equation (1.1) to a fully
nonlinear equation for which the convexity is a natural condition
for the ellipticity of the equation. Let $u$ be a smooth,
uniformly convex function defined in a convex domain $\Om\subset
\R^n$. The Legendre transform of $u$, $u^*$, is a smooth,
uniformly convex function defined in $\Om^*=Du(\Om)$, given by
$$u^*(y)=\sup\{x\cdot y-u(x):\ \ x\in\Om\}. \tag 1.8$$
The supremum is attained at the unique point $x$ such that
$y=Du(x)$, and $u$ can be recovered from $u^*$ by the same
Legendre transform. If $u$ is a convex solution of (1.1), $u^*$
satisfies the fully nonlinear equation
$$\det D^2 u^*
  =\sum (\delta_{ij}-\frac {y_iy_j}{1+|y|^2})F^{ij}[u^*],\tag 1.9$$
where $\delta_{ij}=1$ if $i=j$ and $\delta_{ij}=0$ otherwise, and
$$F^{ij}[u^*]=\frac {\p }{\p r_{ij}}\det\, r
                   \ \ \ \text{at}\ r=D^2 u^*.$$
It is known that for any uniformly convex domain $\Om$ and any
smooth function $\phi$ on $\pom$, (1.9) has a unique convex
solution $u^*$ in $\Om$ satisfying $u^*=\phi$ on $\pom$, see
Theorem 5.2. By Theorem 5.2 we can construct a sequence of convex
solutions $(u^*_k)$ to (1.9), such that $(u_k)$, the Legendre
transform of $(u^*_k)$, converges to a complete convex solution
$u$ of (1.1) and the level set $\{x\in\R^n:\ u(x)=1+\inf u\}$ is
not a sphere.

This paper is arranged as follows. In Section 2 we prove estimate
(1.5) and Theorem 1.3. In Section 3 we prove Theorem 1.1. In
Section 4 we prove the case $\sigma=0$ of Theorem 1.2. The case
$\sigma =1$ of Theorem 1.2 will be proved in Section 5. The final
Section 6 contains some applications. We first prove Corollary
1.1, then discuss implications of Theorem 1.3, and finally mention
a few unsolved problems related to our Theorem 1.1-1.3.

\vskip5pt

\noo{\it Recent developments.} \ A major advance, after the paper
was finished in early 2003, has been made by Huisken and Sinestrari
[29], in which they studied the mean curvature flow with surgeries
of 2-convex hypersurfaces in $\R^{n+1}$ for $n\ge 3$. They proved
that at any point with large curvature, the hypersurface after
normalization must be very close to the cylinder $S^{n-1}\times
\R^1$ or a convex cap. Very recently, in another development, the
author, together with Weimin Sheng [30], found a new proof for the
following result of White [27,28], that is for the mean convex flow
up to the first time singularity, a blow-up sequence converges along
a subsequence to a convex mean curvature flow, and the grim reaper
is not a blow-up solution. This proof is based on the curvature
pinching of Huisken and Sinestrari [14,15].

\vskip5pt

\noo{\it Acknowledgement.} \  The author would like to thank Gerhard
Huisken for bringing the problem to his attention. He also wishes to
thank his colleagues Ben Andrews for helpful discussions, Neil
Trudinger for pointing out Proposition 3.1 to him, and in particular
John Urbas for discussions on the proof of Theorem 5.2.

\vskip30pt

\baselineskip=13.6pt
\parskip=2.5pt

\centerline {\bf 2. Level set estimates and proof of Theorem 1.3}

\vskip10pt

In this section we prove estimate (1.5) and Theorem 1.3. Our proof
of (1.5) involves elementary, but delicate analysis. To illustrate
our idea, we first prove it in the dimension two case, then prove it
for higher dimensions. For clarity we divide this section into 3
subsections. In \S2.1 we prove (1.5) for $n=2$. In \S2.2 we prove
(1.5) for $n>2$. In \S2.3 we prove Theorem 1.3.

Let $u$ be a complete convex solution of (1.2). For any constant
$h>0$ we denote
$$\eqalign
{\Gamma_h & =\Gamma_{h, u}=\{x\in\R^n:\ u(x)=h\},\cr
 \Om_h & =\Om_{h, u}=\{x\in\R^n:\ u(x)<h\} .\cr } \tag 2.1$$
Then $\Om_h\subset \Om_{h+\eps}$ for any $\eps>0$. Let $\kappa$
denote the mean curvature of the level set $\Gamma_h$. We have
$$\align
\Cal L_\sigma [u]
 &=\kappa u_\gamma
    +\frac {\sigma u_{\gamma\gamma}}{\sigma+u_\gamma^2} \tag 2.2 \\
 & \ge \kappa u_\gamma = \Cal L_0 [u] ,\\
     \endalign $$
where $\gamma$ is the unit outward normal to $\Om_{h, u}$, and
$u_{\gamma\gamma}=\gamma_i\gamma_ju_{ij}$.

\vskip5pt

We may assume that $\Om_h$ does not contain a straight line. For if
$\Om_{h_0}$ contains a straight line, then for all $h\ge h_0$, by
convexity we have the splitting $\Om_h=\Om'_h\times\R^1$ for some
convex set $\Om'_h\subset\R^{n-1}$.  The convexity implies that $u$
is a function of $x_1, \cdots, x_{n-1}$ if the straight line is
parallel to the $x_n$-axis. Therefore the problem reduces to a lower
dimension case. Furthermore, the graph of $u$, $\M_u$, does not
contain any line segment. For it it does, the analyticity of $u$
when $\sigma>0$, or the constant ranking of $(D^2 u)$ when
$\sigma=0$ [14,15], implies that $\M_u$ contains a straight line.

We also remark that our proof of (1.5) works for general convex
solutions of (1.2). When $\sigma=0$ and the solution is a limit flow
(blow-up solution) to a mean convex flow, one may also use the
noncollapsing result in [27, 28, 30] to give an alternative proof,
see Remark 2.1 below.

\vskip10pt

\noo{\bf 2.1. Proof of (1.5) for n=2}. Let $u$ be a complete convex
solution of (1.2) satisfying $u(0)=0$. We first prove that if
$\Om_1\cap\{x_1=0\}$ is contained in $\{|x_2|\le \beta\}$ for some
small $\beta>0$, then $\Om_h$ is contained in a strip region for any
$h>0$. We prove the result in three lemmas. In the first one we
assume that $\sigma=0$ and $u$ is symmetric in $x_2$. In the second
one we remove the symmetry assumption. In the third one we remove
the condition $\sigma=0$.

\proclaim{Lemma 2.1} Let $u$ be a complete convex solution of (1.2).
Suppose $n=2$, $\sigma=0$, $u(0)=0$, and $u$ is symmetric in $x_2$,
namely $u(x_1, x_2)=u(x_1, -x_2)$. Suppose there is a sufficiently
small $\beta>0$ such that $u(0, \beta)\ge 1$.  Then $u$ is defined
in a strip region $\{|x_2|<C\}$.
\endproclaim

\noo{\it Proof}. Let $\M^+_u$ denote the graph of $u$ in the
half-space $\{x_2\ge 0\}$, and $D$ the projection of $\M^+_u$ on the
plane $\{x_2=0\}$.  Then $\M^+_u$ can be represented as the graph of
a function $g$ in the form $\M^+_u=\{(x_1, x_2, x_3):\ \ x_2=g(x_1,
x_3)\}$, and $g$ is positive, concave, monotone increasing in $x_3$,
defined in $D$, and vanishes on $\p D$. In the following we also
regard the height parameter $h$ as a variable and use $h$ instead of
$x_3$.

For any $h>0$, we denote $g_h (x_1)=g(x_1, h)$  and $D_h=\{x_1:
(x_1, h)\in D\}$. So $g_h$ is a positive, concave function of one
variable, and $D_h$ is an interval in $\R^1$ containing the origin,
$D_h=(-\underline a_h, \ol a_h)$ (here $\ul a_h$ or $\ol a_h$ might
be equal to infinity). Denote $\ol b_h=g_h(0).$

\vskip5pt

\noo {\it Claim 1:} For any given $h>0$, if $\ol a_h, \underline a_h
\ge \ol b_h$, then $\ol a_h\ol b_h \ge \frac \pi{32}h $.

To prove the claim, we assume $\ol a_h\le \underline a_h$, otherwise
we make the change $x_1\to -x_1$. Denote $U_h=\Om_h\cap\{x_1>0\}$.
When $\sigma=0$, the level set $\Ga_{h}$ is evolving at the velocity
equal to its curvature (with time $t=-h$). By the convexity of $U_h$
and the assumption $\ol a_h, \underline a_h \ge \ol b_h$, we have
$\ol a_s, \underline a_s \ge \frac 12 \ol b_s$ for all $s\in (\frac
12 h, h)$. Hence by the concavity of $g$ we have the gradient
estimate $|\frac {d}{d x_1}g_s(0)|\le 2$ for $s\in (\frac 12 h, h)$.
Note that $\frac {d}{ds}|U_s|_{\Cal H^2}$ is equal to the arc-length
of the set of the unit normals to $\Gamma_s\cap\{x_1>0\}$.  Hence
$\frac {d}{ds}|U_s|_{\Cal H^2}\ge \frac \pi 4 $ for $s\in (\frac 12
h, h)$, which implies $|U_h|_{\Cal H^2}\ge \frac \pi 8h$. Here and
below we use $|E|_{\Cal H^k}$ to denote the $k$-dim Hausdorff
measure of the set $E$. By the convexity of $\Om_h$ and the
assumption $\underline a_h\ge \ol a_h$, one sees that $U_h$ is
contained in $(0, \ol a_h)\times (-2\ol b_h, 2\ol b_h)$. Hence $\ol
a_h\ol b_h\ge \frac 14|U_h|_{\Cal H^2}$. We obtain $\ol a_h\ol
b_h\ge \frac\pi {32} h$.

Claim 1 is also true when $\sigma>0$. Indeed, by equation (2.2) and
the convexity we have $\L_0[u]\le 1$, which means $\Ga_{h}$ is
moving at a velocity greater than or equal to its curvature.
Therefore we also have $\frac {d}{ds}|U_s|_{\Cal H^2}\ge \frac \pi 4
$ for $s\in (\frac 12 h, h)$, and so it also follows $\ol a_h\ol b_h
\ge \frac \pi {32} h $.

In particular, the proof implies Claim 1 holds for convex functions
$u$ satisfying $\L_0[u]\le C_1$ for some positive constant $C_1$.
That is if $\L_0[u]\le C_1$, then $\ol a_h\ol b_h \ge C_2h$, where
$C_2$ depends on $C_1$.

\vskip5pt

\noo {\it Claim 2:} Denote $h_k=2^k$, $\ol a_k=\ol a_{h_k}$, $\ol
b_k=\ol b_{h_k}$, $g_k=g_{h_k}$, and $D_k=D_{h_k}$. Then
$$g_k (0)\le g_{k-1}(0)+ 2^{-k/8}\
     \ \ \ \text{for all}\ \ k\ge 1. \tag 2.3$$
Note that Lemma 2.1 follows from Claim 2 immediately. Indeed, let
$\Om_\infty=\bigcup_{h>0}\Om_h$ be the domain of definition of $u$.
By (2.3), $\ol b_k$ is uniformly bounded. Hence by Claim 1,
$\Om_\infty$ is a convex set containing the whole $x_1$-axis. Hence
$\Om_\infty=I\times \R^1$ for some interval $I$ in the $x_2$ axis.
Estimate (2.3) implies that $I\subset (-2, 2)$ (see (2.5) below).
Hence $\Om_\infty$ must be a strip region.

To prove (2.3) we observe that, since $g$ is positive and concave,
$g_k(0)\le h_k g_0(0)\le 2^k\beta$. Hence we may assume that
$g_{k_0}(0)\le 1$ for some sufficiently large $k_0$. By Claim 1, we
have
$$\ol a_k \ge C_0h_k \tag 2.4$$
for $k\le k_0$ with $C_0=\frac \pi {32}$. We prove (2.3) by
induction, starting at $k=k_0$.

Suppose by induction that (2.3) holds up to $k$. Then by induction,
$$g_k(0)\le g_{k_0}(0)+\sum_{j=k_0}^k 2^{-j/8}\le 2.\tag 2.5$$
By the concavity of $g$ and since $g\ge 0$, we have $g_{k+1}(0)\le
2g_k(0)\le 4$. By Claim 1, $\ol a_{k+1}\ge \frac{\pi}{128} h_k$.
Hence (2.4) holds at $k+1$ with $C_0=\frac\pi{128}$.

To prove (2.3) at $k+1$, denote
$$\align
 L_k &=\{x_1\in\R^1:\ \ -\frac {C_0}{4} h_k<x_1<\frac {C_0}{4} h_k\},\\
  Q_k &=L_k\times [h_k, h_{k+1}]\subset D,\\
  \endalign $$
where $C_0$ is given in (2.4). Since $g>0$, $g$ is concave, and
$g(x_1, h)$ is defined in $2L_k$ for $h\ge \frac 12 h_k$, we have
$$\align
\sup\{g(x_1, h):\ (x_1, h)\in Q_k\}
  & \le 2\sup\{g(x_1, h_k):\ x_1\in L_k\}\tag 2.6\\
  & \le 4 g(0, h_k)=4g_k(0)\le 8.\\
  \endalign $$
Moreover,  for any $(x_1, h)\in Q_k$,
$$|\p_h g(x_1, h)|
 \le \frac {g(x_1, h)-g(x_1, h_{k-1})}{h-h_{k-1}}
  < \frac {g(x_1, h)}{h_k-h_{k-1}}
   \le \frac {16}{h_k}.\tag 2.7$$
Similarly,
$$|\p_{x_1} g(x_1, h)|\le \frac {g(x_1, h)}{\ol a_h-|x_1|}
          \le \frac {2g(x_1, h)}{h_k}\le \frac{16}{h_k}
 \ \ \ \forall\ (x_1, h)\in Q_k.\tag 2.8$$
From the above gradient estimates and the concavity of $g$, the
average in $Q_k$ of the second derivative satisfies
$$|\p^2_{x_1} g(x_1, h)|
  \approx  \sup_{x_1\in L_k} {|\p_{x_1} g(x_1, h)|}/|L_k|_{\Cal H^1}
                  \approx  h_k^{-2}.\tag *$$
Here $a\approx b$ means $a\le Cb$ and $b\le Ca$ for some constant
$C$. This simple observation is {\it critical} for the proof of
(2.3), and actually the following weaker version is sufficient,
$$|\p^2_{x_1} g (x_1, h)|
  \le Ch_k^{-5/4}\ \ \ (x_1, h)\in Q_k\tag 2.9$$
for some fixed constant $C$. Indeed, by equation (2.2), we have
$\kappa u_\gamma = 1$ (when $\sigma=0$). Note that $u_\gamma\approx
(\p_h g)^{-1}$ and $\kappa\approx \p_{x_1x_1}g$ when $(x_1, h)\in
Q_k$. Hence if (2.9) holds, we have
$$|\p_h g(0, h)|\le Ch_k^{-5/4}
              \ \ \text{for}\ \ h\in (h_k, h_{k+1}). \tag 2.10$$
It follows that
$$g_{h_{k+1}}(0)-g_{h_k}(0)=g(0, h_{k+1})-g(0, h_k)
      \le Ch_k^{-1/4}\le h_k^{-1/8}\tag 2.11$$
when $k$ is large. We obtain (2.3).

However we have not proved the estimate (2.9) pointwise, even in the
special case $\sigma=0$. But we observe that the set where $g$ does
not satisfies (2.9) is very small, and the concavity of $g$ ensures
that this small set does not harm the estimate (2.3).

Denote
$$\chi = \{(x_1, h)\in Q_k:\
   |\p_{x_1}^2 g(x_1, h)| \ge h_k^{-5/4}\} .$$
If $\chi$ is empty, (2.3) is proved in (2.11) above. If
$\chi\ne\emptyset$, we proceed as follows (the argument also applies
to the case $\chi$ is empty). For any $h\in (h_k, h_{k+1})$, by the
gradient estimates and the concavity of $g$, an integration by parts
gives
$$|\{x_1\in L_k:\ (x_1, h)\in\chi\}|_{\Cal H^1} h_k^{-5/4}
   \le \big{|}\int_{L_k} \p_{x_1x_1} g\big{|}
   \le 2 \ {\sup}_{L_k} |\p_{x_1} g|
   \le Ch_k^{-1}.$$
Taking integration from $h=h_k$ to $h=h_{k+1}$ we obtain
$|\chi|_{\Cal H^2} h_k^{-5/4}\le C$, namely
$$|\chi|_{\Cal H^2}\le C h_k^{5/4} . \tag 2.12$$
We say $\chi$ is a small set as the ratio $|\chi|_{\Cal
H^2}/|Q_k|_{\Cal H^2}=O(h_k^{-3/4})$ is small.

For any given $y_1\in L_k$ we denote $\chi_{y_1}=
\chi\cap\{x_1=y_1\}$. From (2.12) and by the Fubini Theorem, there
is a set $\wtt L\subset L_k$ with measure $|\wtt L|_{\Cal H^1} <
h_k^{1/2}$ such that for any $y_1\in  L_k-\wtt L$,
$$|\chi_{y_1}|_{\Cal H^1}\le C h_k^{3/4}. \tag 2.13$$
For any $y_1\in L_k-\wtt L$, we have
$$
g(y_1, h_{k+1})-g(y_1, h_k)
 =\int_{\chi_{y_1}}\p_h g(y_1, h)dh
  +\big(\int_{h_k}^{h_{k+1}}-\int_{\chi_{y_1}}\big)\p_h g(y_1, h)dh.$$
By (2.13) and the gradient estimate (2.7), the first integral is
bounded by $Ch_k^{-1/4}$. For the second one, similarly as in (2.10)
we have $\p_h g\le Ch_k^{-5/4}$ for any point $(x_1, h)\not\in\chi$.
Hence the second integral is also bounded by $Ch_k^{-1/4}$.
Therefore we obtain
$$g(y_1, h_{k+1})-g(y_1, h_k)\le Ch_k^{-1/4}. \tag 2.15$$

From (2.15) we get (2.3) immediately. Indeed, replacing $x_1$ by
$-x_1$ if necessary, we assume that $\p_{x_1}g(0, h_k)\le 0$, so
that $g(x_1, h_k)\le g(0, h_k)$ for all $x_1\ge 0$. Since $|\wtt
L|<h_k^{1/2}$, the set $[0, h_k^{1/2}]-\wtt L$ is not empty. Let
$y_1\in [0, h_k^{1/2}]-\wtt L$. By (2.15) we obtain
$$g(y_1, h_{k+1})
       \le g(y_1, h_k)+Ch_k^{-1/4}
              \le g(0, h_k)+Ch_k^{-1/4}.\tag 2.16 $$
Note that $g$ is positive and concave, and that $g$ is defined on
the interval $[0, \ol a_{k+1}]$ with $\ol a_{k+1}\ge Ch_{k+1}$
($C=\frac{\pi}{128}$ as established before). We have
$$\frac {g_{k+1}(0)}{g_{k+1}(y_1)}
 \le \frac {\ol a_{k+1}}{\ol a_{k+1}-|y_1|}=1+Ch_{k+1}^{-1/2}$$
By (2.16) it follows (recall our notation, $g(y_1,
h_{k+1})=g_{k+1}(y_1)$),
$$g_{k+1}(0)\le g_k(0)+Ch_k^{-1/4}. $$
We obtain (2.3). Lemma 2.1 is proved. \hfill$\square$

\vskip10pt

In Lemma 2.1, we don't need to assume that the sub-level set $\Om_h$
is compact, nor we assume the whole sub-level set $\Om_{h\,|\, h=1}$
is contained in a strip region. The assumption that $u(0, \beta)\ge
1$ for sufficiently small $\beta>0$ implies that
$\Om_1\cap\{x_1>0\}$ if $u_{x_1}(0, \beta)>0$, or
$\Om_1\cap\{x_1<0\}$ if $u_{x_1}(0, \beta)<0$, is contained in a
strip $\{|x_2|\le \beta\}$. Next we remove the assumption that $u$
is symmetric in $x_2$.

\proclaim{Lemma 2.2} Let $u$ be a complete convex solution of (1.2).
Suppose $n=2$, $\sigma=0$, $u(0)=0$, and there is a sufficiently
small $\beta>0$ such that $u(0, \beta)\ge 1$ and $u(0, -\beta)\ge
1$. Then $u$ is defined in a strip region.
\endproclaim

Note that the assumption $u(0, \beta)\ge 1$ and $u(0, -\beta)\ge 1$
is equivalent to that $\Om_1\cap\{x_1=0\}\subset\{|x_2|<\beta\}$.
Before proving Lemma 2.2, we state a property of convex domain, due
to F. John [21], which is frequently used in the study of convex
bodies and Monge-Amp\`ere equations.

\proclaim{Proposition 2.1} Let $\Om$ be a bounded, convex domain in
$\R^n$, $n\ge 2$. Then among all (solid) ellipsoids containing
$\Om$, there is a unique ellipsoid $E$ of smallest volume such that
$$\frac 1n E\subset \Om\subset E,\tag 2.17$$
where $\alpha E$ is the $\alpha$-dilation of $E$ with respect to its
center.
\endproclaim

We call $E$ the {\it minimum ellipsoid} of $\Om$ (it is a (solid)
ellipse when $n=2$). By a rotation of the coordinates, we may assume
that $E$ is given by $E=\big\{{\sum}_{i=1}^n\big(\frac {x_i-x_{0,
i}}{r_i} \big)^2<1\big\}$, where $x_0=(x_{0, 1}, \cdots, x_{0, n})$
is the center of $E$. We can make the linear transform
$y_i=(x_i-x_{0,i})/r_i+x_{0,i}$, $i=1, \cdots, n$, such that $E$
becomes the unit ball $B_1(x_0)$ and $B_{1/n}(x_0)\subset
T(\Om)\subset B_1(x_0)$.

\vskip10pt

\noo{\it Proof of Lemma 2.2}. Let $R=10^3$ and let $E$ be the
minimum ellipsoid of $\Om_1\cap B_R(0)$. By a rotation of
coordinates we assume the axial directions of $E$ coincide with
those of the coordinate system (we don't need to assume the center
of $E$ is at the origin).

The proof is similar to that of Lemma 2.1. We indicate the necessary
changes. Let $\M_u$ be the graph of $u$, which consists of two
parts, $\M_u=\M^+\cup\M^-$, where $\M^+=\{(x, u(x))\in\R^3:\
\p_{x_2} u(x) \ge 0\}$ and $\M^-=\{(x, u(x))\in\R^3:\ \p_{x_2} u(x)
\le  0\}$. Then $\M^\pm$ can be represented as graphs of functions
$g^\pm$ in the form $x_2=g^\pm(x_1, x_3)$,  $(x_1, x_3)\in D$ and
$D$ is the projection of $\M_u$ on the plane $\{x_2=0\}$. The
functions $g^+$ and $g^-$ are respectively concave and convex, and
we have $x_3=u(x_1, g^\pm(x_1, x_3))$. Denote
$$g=g^+-g^-.\tag 2.18$$
Then $g$ is a positive, concave function in $D$, vanishing on $\p
D$. For any $h>0$ we also denote $g_h (x_1)=g(x_1, h)$,
$g^\pm_h(x_1)=g^\pm(x_1, h)$, and $D_h=\{x_1\in\R^1: (x_1, h)\in
D\}$. Then $g_h$ is a positive, concave function in $D_h$, vanishing
on $\p D_h$, and $D_h=(-\ul a_h, \ol a_h)$ is an interval containing
the origin. As before we denote $\ol b_h=g_h(0)$.

\vskip5pt

\noo {\it Claim 1:} Suppose $\ol a_h, \underline a_h \ge \ol b_h$.
Then $\ol a_h\ol b_h \ge \frac\pi{32}h $.

The claim can be proved in the same way as in Lemma 2.1, by
observing that the gradient estimate $|\frac {d}{d x_1}g_s(0)|\le 2$
also implies that the arc-length of the set of the unit normals to
$\Gamma_s\cap\{x_1>0\}$ is greater than $\frac \pi 4$. Hence we also
have $\frac {d}{ds}|U_s|_{\Cal H^2}\ge \frac \pi 4 $ for $s\in
(\frac 12 h, h)$.

\noo {\it Claim 2:} Denote $h_k=2^k$, $\ol a_k=\ol a_{h_k}$, $\ol
b_k=\ol b_{h_k}$, $g_k=g_{h_k}$, and $D_k=D_{h_k}$. We have
$$g_k (0)\le g_{k-1}(0)+ 2^{-k/8}\
     \ \ \ \text{for all}\ \ k\ge 1. \tag 2.19$$
Lemma 2.2 follows from Claim 2 immediately. Indeed, let $P$ be the
projection of the graph $\M_g$ on the plane $\{x_3=0\}$. Then $P$ is
convex set containing the $x_1$-axis. Hence $P=I\times \R^1$ for
some interval $I$. Estimate (2.19) implies that $g(0, h)\le 2$ for
all $h$ (due to (2.5)), so we have $I\subset [0, 2]$. Hence $\M_g$
is contained in the strip $\{(x_1, x_2, x_3)\in\R^3:\ 0\le x_2\le
2\}$. By (2.18), $\M_u$ is also contained in a strip region
$\{|a_1x_1+a_2x_2|<2\}$, where $(a_1, a_2, 0)$ is a unit vector in
$\R^3$ with $a_1$ small and $a_2$ close to 1. We can make $a_1$ as
small as we want, provided the constant $R$ at the beginning of the
proof is sufficiently large.

The proof of (2.19) is similar to that of (2.3). But the argument
from (2.9) to (2.10) needs to use the equation (2.2). Therefore we
need to consider $g^+$ and $g^-$ instead of $g$.

First we establish (2.4)-(2.8) in the same way as in Lemma 2.1. Let
$L_k, Q_k$ and $\chi$ be as in Lemma 2.1. We also denote
$$\chi^\pm = \{(x_1, h)\in Q_k:\
   |\p_{x_1}^2 g^\pm (x_1, h)| \ge h_k^{-5/4}\} .$$
Then both $\chi^+$ and $\chi^-$ are subsets of $\chi$. For any $h\in
(h_k, h_{k+1})$, by (2.8) and recalling that $L_k=(-\frac14 C_0 h_k,
\frac 14 C_0 h_k)$, we have
$$\align
|\{x_1\in L_k:\ (x_1, h)\in\chi^+\}|_{\Cal H^1} h_k^{-5/4}
   & \le \big{|}\int_{L_k} \p_{x_1x_1} g^+\big{|}\\
   &= |\p_{x_1} g^+(\frac 14 C_0h_k, h)-\p_{x_1} g^+(-\frac 14 C_0h_k, h)|\\
   & \le 2 \ {\sup}_{L_k} |\p_{x_1} g|
   \le Ch_k^{-1},\\
   \endalign $$
where the second inequality is due to that $g=g^+-g^-$, $g^+$ is
concave and $g^-$ is convex. Hence $|\chi^+|_{\Cal H^2}\le C
h_k^{5/4}$. Similarly we have $|\chi^-|_{\Cal H^2}\le C h_k^{5/4}$.

For any given $y_1\in L_k$, denote $\chi^\pm_{y_1}= \chi^\pm \cap
\{x_1=y_1\}$. Then there is a set $\wtt L^\pm\subset L_k$ with
measure $|\wtt L^\pm|_{\Cal H^1} < h_k^{1/2}$ such that for any
$y_1\in L_k-\wtt L^\pm$, we have $|\chi^\pm_{y_1}|_{\Cal H^1}\le C
h_k^{3/4}$.

For any given $y_1\in L_k-(\wtt L^+\cup\wtt L^-)$, we have
$$g(y_1, h_{k+1})  -g(y_1, h_k)
 =g^+(y_1, h_{k+1})  -g^+(y_1, h_k)
     +|g^-(y_1, h_{k+1})  -g^-(y_1, h_k)|. $$
In the following we estimate $g^+(y_1, h_{k+1})  -g^+(y_1, h_k)$.
The estimate also applies to $|g^-(y_1, h_{k+1}) -g^-(y_1, h_k)|$.
We have
$$\align
 & g^+(y_1, h_{k+1})  -g^+(y_1, h_k)
 = \int_{h_k}^{h_{k+1}} \p_h g^+(y_1, h) dh \tag 2.20\\
 & =\int_{\chi^+_{y_1}}\p_h g^+(y_1, h)dh
  +\big(\int_{h_k}^{h_{k+1}}
     -\int_{\chi^+_{y_1}}\big)\p_h g^+(y_1, h)dh.\\
     \endalign $$

For the first integral on the right hand side, note that
$g=g^+-g^-$, $g^+(y_1, h)$ is concave and increasing in $h$, and
$g^-$ is convex and decreasing in $h$. Hence $\p_h g^+\le \p_h g$.
By the gradient estimate (2.7) and recalling that
$|\chi^+_{y_1}|_{\Cal H^1}\le C h_k^{3/4}$, we have
$$\int_{\chi^+_{y_1}}\p_h g^+(y_1, h)dh\le Ch_k^{-1/4}.$$
To estimate the second integration, we first introduce a mapping
$\Cal T: p\to q$ as follows. For a point $p=(x_1, h)\in D$, there is
a corresponding point $P=(x_1, x_2, h)$ on the level set $\Gamma_h$,
where $x_2=g^+(x_1, h)$, such that $p$ is the projection of $P$ on
the plane $\{x_2=0\}$.  Let $q=(x_1, x_2)$ be the projection of $P$
on the plane $\{x_3=0\}$.

By equation (2.2), we have $\kappa u_\gamma = 1$. Note that when
$p=(x_1, h)\in Q_k$, the normal $\gamma$ of the level set
$\Gamma_h\subset \R^2$ at the point $q=\Cal T(p)\in \Gamma_h$
satisfies
$$|\gamma - e_2|<\eps\tag 2.21$$
for some small constant $\eps>0$, where $e_2=(0, 1)$. This is
because by induction, $\Om_{h_k}\cap B_R$ is contained in a strip
region (see discussion after (2.19)) and the axial directions of the
minimum ellipsoid of $\Om_{h_k}\cap B_R$ is a small perturbation of
the axial directions of the coordinates, where $R$ is the constant
introduced at the beginning of the proof.  Therefore we have
$$\cases
 (\p_h g^+)^{-1} =(1+\eps_1) u_\gamma,\\
 \p_{x_1x_1}g^+ =(1+\eps_2)\kappa,\\
 \endcases\tag 2.22$$
where $\eps_1, \eps_2$ are small constants provided $R$ is
sufficiently large. Hence
$$|\p_h g^+(y_1, h)|\le C|\p_{x_1x_1}g^+|\le Ch_k^{-5/4}
 \ \ \forall\ (y_1, h)\not\in \chi^+.\tag 2.23 $$
It follows that
$$\big(\int_{h_k}^{h_{k+1}}
     -\int_{\chi^+_{y_1}}\big)\p_h g^+(y_1, h)dh \le Ch_k^{-1/4}.$$
Combining the above two estimates we obtain
$$g^+(y_1, h_{k+1})-g^+(y_1, h_k)\le Ch_k^{-1/4}. $$
Similarly we have $|g^-(y_1, h_{k+1})-g^-(y_1, h_k)|\le
Ch_k^{-1/4}$. Therefore we obtain (2.19) just as we prove (2.3) from
(2.15). \hfill$\square$

\vskip5pt

Next we remove the condition $\sigma=0$ in Lemma 2.2.

\proclaim{Lemma 2.3} Let $u$ be a complete convex solution of (1.2).
Suppose $n=2$, $u(0)=0$, and there is a sufficiently small $\beta>0$
such that $u(0, \beta)\ge 1$ and $u(0, -\beta)\ge 1$. Then $u$ is
defined in a strip region.
\endproclaim

\noo{\it Proof}. We can follow the proof of Lemma 2.2 until (2.20)
without any change. The estimate for the second integral on the
right hand side of (2.20) used the equation $\kappa u_\gamma =1$.
But when $\sigma\ne 0$, equation (2.2) contains an extra term
$\frac{\sigma u_{\gamma\gamma}}{\sigma+u_\gamma^2}$. To handle this
extra term, we need to divide the integral (2.20) into three parts,
$$\align
g^+(y_1, h_{k+1})-g^+(y_1, h_k)
  & =\int_I \p_h g^+(y_1, h)dh\tag 2.24\\
  & =\left(\int_{I_1}+\int_{I_2}+\int_{I_3}\right)\p_h g^+(y_1, h)dh\\
  \endalign$$
where $I=(h_k, h_{k+1})$,
$$\align
 I_1 & =\chi^+_{y_1},\\
 I_2 & = \{h\in I-I_1:
  \ \frac{\sigma u_{\gamma\gamma}(q)}{\sigma+u_\gamma^2(q)} \le \frac 12\},\\
 I_3 & =I-(I_1\cup I_2), \\
 \endalign $$
where $q=\Cal T(p)$, $p=(y_1, h)$, and $\Cal T$ is the mapping
introduced after (2.20).

Similarly as in Lemma 2.2, we have $\p_h g^+\le \p_h g\le C/h$ for
$h\in (h_k, h_{k+1})$ and the first integral $\int_{I_1}\p_h
g^+(y_1, h)dh\le Ch_k^{-1/4}$.

For the second one, noting that when $\frac{\sigma
u_{\gamma\gamma}}{\sigma+u_\gamma^2} <\frac 12$, we have, by (2.22),
$$(\p_h g^+)^{-1} \p_{x_1x_1}g^+\approx \kappa u_\gamma \ge \frac 12.$$
Hence $\p_h g^+\le C \p_{x_1x_1}g^+$ and we obtain $\int_{I_2}\p_h
g^+(y_1, h)dh\le Ch_k^{-1/4}.$

To estimate the third integral in (2.24), note that for any point
$p=(y_1, h)$ with $h\in I_3$, we have
$$\cases
u_\gamma(q)  =u_{x_2}(q)(1+\eps_1),\\
u_{\gamma\gamma}(q)  =u_{x_2x_2}(q)(1+\eps_2)+o(u_{x_2})\\
 \endcases \tag 2.25$$
at the point $q=\Cal T(p)$, for some small constants $\eps_1$ and
$\eps_2$. The first formula is due to (2.21). To verify the second
one in (2.25), one chooses a coordinate system $(z_1, z_2)$ such
that $q$ is the origin, $\gamma$ is in the $z_2$-axis and $\Gamma_h$
is locally given by $z_2=\eta(z_1)$. One then differentiates $u(z_1,
\eta(z_1))=h$ twice to obtain $u_{z_1z_1}+u_{z_2}\kappa=0$. Recall
that when $p=(y_1, h)$ with $h\in I_3$, $\kappa\le C h_k^{-5/4}$.
Hence $u_{z_1z_1}=o(u_{z_2})$, from which one easily obtains (2.25).

Since for any $p\in \{y_1\}\times I_3$,
 $\frac{\sigma u_{\gamma\gamma}(q)}{\sigma+u_\gamma^2(q)} \ge \frac 12$.
Therefore by (2.25) we have
$$\frac {\sigma u_{x_2x_2}(q)}{\sigma+u_{x_2}^2(q)}
                        \ge \frac 13.$$
Notice that $0\le\sigma\le 1$. Hence
$$u_{x_2x_2}\ge \frac 14(\sigma+u_{x_2}^2)\ge \frac 14 u_{x_2}^2. \tag 2.26$$

Now by the relation $h=u(y_1, x_2)$ and $x_2=g^+(y_1, h)$, we have
$h=u(y_1, g^+(y_1, h))$. As $y_1$ is fixed, we can regard $u$ and
$g^+$ as functions of one variable. Differentiating in $h$ gives
$1=u'(g^+)'$, differentiating twice we get
$0=u''{(g^+)'}^2+u'(g^+)''$. Hence
$$(g^+)''=-\frac{u''}{u'} {(g^+)'}^2=-u''{(g^+)'}^3.$$
By (2.26) we then obtain
$$(g^+)''\le -\phi(h) (g^+)',\tag 2.27$$
where $\phi(h)=\frac 14$ if $h\in I_3$ and $\phi(h)=0$ otherwise.
Observing that $(g^+)'>0$, we obtain
$$\int_{h_k}^h\frac {(g^+)''}{(g^+)'}
       \le -\int_{h_k}^h\phi
             =-\frac 14|I_{3, h}|\ \ \ \forall\ h\in (h_k, h_{k+1}), $$
namely
$$\log (g^+)'(h)\le \log (g^+)'(h_k)-\frac 14 |I_{3,h}|,$$
or equivalently
$$(g^+)'(h)\le (g^+)'(h_k)e^{-|I_{3,h}|/4},\tag 2.28$$
where $I_{3,h}=I_3\cap[h_k, h]$. Since $I_3=\cup_k (a_k, b_k)$ is
the union of intervals, and $g^+$ is increasing, the third integral
in (2.24) is equal to $\text{osc}_{I_3} g^+=\sum_k
g^+(b_k)-g^+(a_k)$. We have
$$\align
\text{osc}_{I_3} g^+
 &\le (g^+)'(h_k)\int_{I_3}e^{-|I_{3,h}|/4} \tag 2.29\\
 &\le (g^+)'(h_k)\int_{h_k}^{h_{k+1}}e^{-(h-h_k)/4}\\
 &\le 2(g^+)'(h_k)\le \frac {C}{h_k}.\\
 \endalign  $$
This completes the proof.\hfill$\square$

\vskip5pt

Next we prove an auxiliary lemma.

\proclaim{Lemma 2.4} Let $u$ be a complete convex solution of (1.2).
Suppose $n=2$, $u(0)=0$,  $\delta:=\inf\{|x|:\ x\in\Ga_1\}$ is
attained at $x_0=(0, -\delta)\in\Gamma_1$, and $\delta>0$ is
sufficiently small. Then $D_1$ contains the interval $(-R, R)$ with
$$R\ge (-\log\delta-C)^{1/2},\tag 2.30$$
where $C>0$ is independent of $\delta$, $D_h$ is the set introduced
in the proof of Lemma 2.2.
\endproclaim

\noo{\it Proof}. Suppose near $x_0$, $\Gamma_1$ is given by
$$x_2=g(x_1).$$
Then $g$ is a convex function, $g(0)=-\delta$, and $g'(0)=0$. Let
$a, b>0$ be two constants such that $g(a)=0$ and $g'(b)=1$. To prove
(2.30) it suffices to prove
$$b\ge (-\log\delta-C)^{1/2} .\tag 2.31$$

For any $y=(y_1, y_2)\in \Gamma_1$, where $y_1\in [0, b]$, let $\xi
= y/|y|$. By the convexity of $u$,
$$u_\xi(y) \ge \frac {u(y)-u(0)}{|y|}= \frac {1}{|y|}.$$
Let $\th$ denote the angle between $\xi$ and the tangential vector
$\frac {1}{\sqrt{1+{g'}^2}} (1, g')$ of $\Ga_1$ at $y$. Then
$$\align
\cos\th & =\frac {\xi_1+\xi_2g'(y_1)}{\sqrt{1+{g'}^2}} , \\
 \sin\th & =\sqrt{1-\cos^2\th}
           = \frac {\xi_1g'-\xi_2}{\sqrt{1+{g'}^2}}.\\
 \endalign $$
Hence
$$u_\gamma (y)= {u_\xi(y)}/{\sin\th}
            \ge \frac {\sqrt{1+{g'}^2}}{y_1g'-y_2}, $$
where $\gamma$ is the normal of the sub-level set $\Om_1=\{u<1\}$.
By $\Cal L_0[u]\le 1$, we obtain,
$$\frac {g''}{(1+{g'}^2)^{3/2}}
     \frac {\sqrt{1+{g'}^2}}{y_1 g'-y_2}
     \le \kappa u_\gamma(y) \le 1,$$
where $\kappa$ is the curvature of the level set $\Ga_1=\{u=1\}$.
Hence
$$\align
g''(y_1)
& \le (1+{g'}^2)(y_1g'-y_2)\tag 2.32\\
& \le \cases
2(y_1g'+\delta)\ \ \text{if}\ \ y_2\le 0,\\
2y_1g' \ \ \ \text{if}\ \ y_2\ge 0,\\
\endcases\\
\endalign$$
where $y_2=g(y_1)$ and $g'(y_1)\le 1$ for $y_1\in (0, b)$. We
consider the equation
$$\rho''(t) = \cases
2(t\rho'(t)+\delta)\ \ \text{if}\ \ \rho(t)\le 0\\
2t\rho'(t) \ \ \ \text{if}\ \ \rho(t)\ge 0\\
\endcases$$
with the initial condition $\rho(0)=-\delta$ and $\rho'(0)=0$. Let
$\alpha>0$ such that $\rho(\alpha)=0$. Then for $t\in (0, \alpha)$
we have
$$\rho'(t)=2\delta e^{t^2}\int_0^t e^{-s^2}ds. $$
Hence we have $C_1\le \alpha\le C_2$ and $\rho'(\alpha)\le
C_2\delta$ for some constants $C_1, C_2$. Let $\beta>\alpha$ such
that $\rho'(\beta)=1$. Consider the equation
$$\rho''= 2t\rho'$$
in the interval $(\alpha, \beta)$. Then
$$\log \rho'\big{|}^\beta_\alpha= t^2\big{|}^\beta_\alpha$$
We obtain
$$\beta^2\ge |\log\delta|-C. $$
By the comparison principle we have $g\le \rho$. Hence (2.31) holds.
\hfill$\square$

\proclaim{Theorem 2.1} Let $u$ be an entire convex solution of (1.2)
in $\R^2$. Then
$$u(x)\le C (1+|x|^2), \tag 2.33$$
where the constant $C$ depends only on the upper bound for $u(0)$
and $|Du(0)|$.
\endproclaim

\noo{\it Proof.} By adding a constant to $u$ we may suppose
$u(0)=0$. To prove (2.33) it suffices to prove that $\dist(0,
\Ga_h)\ge Ch^{1/2}$ for all large $h$. By the rescaling
$u_h(x)=\frac 1h u(h^{1/2}x)$ it suffices to prove $\dist(0, \Ga_{1,
u_h})\ge C$. Note that $|Du_h(0)|=h^{-1/2} |Du(0)|\le |Du(0)|$.
Hence by convexity, $\inf_{B_1(0)} u_h$ is uniformly bounded from
below. Note also that $u_h$ satisfies equation (1.2) with $\sigma\to
0$ as $h\to\infty$.

Denote $\delta=:\inf\{|x|:\ x\in\Ga_{1, u_h}\}$. Suppose the infimum
is attained at $x_0=(0, -\delta)$. If $\delta>0$ is sufficiently
small, by Lemma 2.4, $D_1=D_{1, u_h}$ contains the interval $(-R,
R)$, where $R=(-\log\delta-C)^{1/2}$. Let $\delta^*>0$ such that
$u_h(0, \delta^*)=1$. Then $\delta^*$ must also be very small, for
otherwise by convexity the ellipse
$$E=\{(x_1, x_2)\in\R^2:\ \frac {x_1^2}{(R/4)^2}
   +\frac {|x_2-(\delta^*-\delta)/2|^2}
   {[(\delta^*+\delta)/8]^2}<1\}$$
is contained in sub-level set $\Om_{1, u_h}$.

When $\sigma=0$, the level set $\Gamma_{-t, u_h}$ is a solution to
the curve shortening flow (for time $t$ starting at $-1$). Let
$E_{-t}$ be the solution to the curve shortening flow with initial
condition $E_{-1}=E$, where $E$ is the ellipse given above.
Therefore we have the inclusion $E_{-t}\subset \Om_{-t, u_h}$ for
all $t>-1$. It takes the time $T=|E|_{\Cal H^2}$ for the solution
$E_{-t}$ to shrink to a point. Hence we have $\inf_{B_1(0)} u_h\le
1-T$. But when $\delta$ is small and $\delta^*$ has a positive lower
bound (independent of $\delta$), $T=|E|_{\Cal H^2}$ becomes
sufficiently large, which contradicts with the assertion that
$\inf_{B_1(0)} u_h$ is uniformly bounded from below.

When $\sigma\in (0, 1]$, $u_h$ is a solution of (1.2) with
$\sigma\le 1/h$. If there is a sequence $h_k\to\infty$ and
$\delta^*_k\ge \delta^*$ for some $\delta^*>0$ such that $u_{h_k}(0,
\delta^*_k)=1$, we define $E$ as above. Now let $v_\sigma$ be the
solution of $\Cal L_\sigma(v)=1$ in $E$ and $v=1$ on $\p E$. Then
for any given $\delta, \delta^*>0$ and $R>1$, the solution
$v_\sigma$ converges to $v_0$, the solution to $\L_0(v_0)=1$ in $E$
and $v_0=1$ on $\p E$. The level set of $v_0$ is a solution to the
curve shortening flow. Hence $\inf v_\sigma \to-\infty$ as $\delta,
\sigma\to 0$. We also reach a contradiction. \hfill$\square$

\vskip5pt

\proclaim{Corollary 2.1} Let $u$ be a complete convex solution of
(1.2). Then $u$ is either an entire solution, or is defined in a
strip region. In particular, there is no complete convex solution of
(1.2) defined in a half space.
\endproclaim

\noo{\it Proof.}\ \ Let us assume $u(0)=0$. If $u$ is not an entire
solution, then for any $M>1$, there exists $x_0\in\R^n$ such that
$u(x_0)>M|x_0|^2$. Let $u_h(x)=h^{-1}u(h^{1/2}x)$, where $h=u(x_0)$.
Then the distance from the origin to the level set $\Ga_1=\{u_h=1\}$
is less than $M^{-1}$. The proof of Theorem 2.1 then implies that
$u$ is defined in a strip region. \hfill$\square$

\vskip5pt

Note that in the above proof we have used the following lemma.

\proclaim{Lemma 2.5} Let $u_k$ be a sequence of convex solutions of
(1.2) with $\sigma=\sigma_k\in [0, 1]$. Suppose $\sigma_k\to\sigma$
and $u_k\to u$. Then $u$ is a convex solution of (1.2).
\endproclaim

\noo{\it Proof}. Lemma 2.5 is well known if $\sigma_k\equiv 0$ or
$\sigma_k\equiv 1$. If $\sigma_k\to\sigma>0$,  replacing $u_k$ by
$\frac 1{\sigma_k} u(\sqrt{\sigma_k} x)$ we may suppose
$\sigma_k\equiv 1$. We need only to consider the case when
$\sigma_k\to 0$.

In this case we show that $u$ is a viscosity solution of $\Cal
L_0[u]=0$. Indeed, since $\L_{\sigma_k}[u_k]=1$, by convexity we
have $\L_0[u_k]\le 1$ and so $\L_0[u]\le 1$. On the other hand, for
any fixed $\hat\sigma>0$, by convexity we have
$\L_{\hat\sigma}[u_k]\ge \L_{\sigma_k}[u_k]=1$ if $k$ is
sufficiently large such that $\sigma_k<\hat\sigma$. Hence
$\L_{\hat\sigma}[u]\ge 1$. As $\hat\sigma>0$ is arbitrary, we have
$\L_0[u]\ge 1$. \hfill$\square$

\vskip5pt

\noo{\bf Remark 2.1}. When $\sigma=0$ and $u$ is a blow-up solution
(limit flow) to a given mean convex flow, by a compactness argument,
together with Lemma 2.4 and the proof of Theorem 2.1, one sees that
(1.5) also follows from the non-collapsing in [27,28,30]. That is if
(1.5) is not true, there exists a sequence of blow-up solutions
$u_k$ to a given mean convex flow such that
$w_k(x):=k^{-1}u_k(k^{1/2}x)$ converges to a multiplicity two plane.
But a multiplicity two plane does not occur as a blow-up solution
[27, 28, 30].

\vskip10pt

\noo{\bf 2.2. Proof of (1.5) for n$>$2}. In this subsection we
extend the results in \S2.1 to high dimensions.

Let $u$ be a complete convex solution of (1.2). Let $\M_u$ denote
the graph of $u$, and $D$ the projection of $\M_u$ on the plane
$\{x_n=0\}$. We divide $\M_u$ into two parts, $\M_u=\M^+\cup\M^-$,
where $\M^\pm=\{(x, u(x))\in\R^{n+1}:\ \p_{x_n} u(x) \gtreqless
0\}$. Then $\M^+$ and $\M^-$ can be represented respectively as
graphs of the form $x_n=g^+(x', x_{n+1})$ and $x_n=g^-(x',
x_{n+1})$, where $x'=(x_1, \cdots, x_{n-1})$, $(x', x_{n+1})\in D$.
The functions $g^+$ and $g^-$ are respectively concave and convex,
and satisfy the relation $x_{n+1}=u(x', g^\pm(x', x_{n+1}))$. As
before we denote $g=g^+-g^-$. Then $g$ is a positive, concave
function in $D$, vanishing on $\p D$.

For any $h>0$ we also denote $g_h (x')=g(x', h)$,
$g^\pm_h(x')=g^\pm(x', h)$, and $D_h=\{x'\in\R^{n-1}: (x', h)\in
D\}$. Then $g_h$ is a positive, concave function in $D_h$, vanishing
on $\p D_h$, and $D_h$ is a convex domain in $\R^{n-1}$ containing
the origin. Hence $\p D_h$ can be represented as a radial graph of a
positive function $a_h$ on $S^{n-2}$, $\p D_h=\{p\cdot a_h(p):\ p\in
S^{n-2}\}$, where
$$a_h(p)=\sup\{t:\ \ tp\in D_h\}, \ \ \ p\in S^{n-2}.$$
Denote
$$\align
\ol a_h & =\inf\{a_h(p):\ \ p\in S^{n-2}\},\\
\ol b_h &= g_h(0).\\
\endalign $$
We want to extend Lemma 2.3 to high dimensions, that is if $\ol
b_1(0)$ is small, then $u$ is defined in a strip region. First we
prove a lemma which corresponds to Claim 1 in the proof of Lemmas
2.1 and 2.2.

\proclaim{Lemma 2.6} Let $u$ be a complete convex solution of (1.2)
satisfying $u(0)=0$. Suppose $\ol a_h \ge \ol b_h$. Then there is a
positive constant $C_n$, depending only on $n$, such that
$$\ol a_h\ol b_h \ge C_nh . \tag 2.34$$
\endproclaim

\noo{\it Proof.} When $n=2$, (2.34) was proved in Claim 1 in Lemmas
2.1 and 2.2. When $n\ge 3$, we reduce (2.34) to the case $n=2$.

Assume that $\ol a_h=a_h(p)$ for $p= (1, 0, \cdots, 0)$. Observing
that $\ol a_h\ol b_h$ is propositional to the area of the section
$\{x\in\Om_h: x_1>0, x_2=\cdots=x_{n-1}=0\}$, we can prove (2.34) by
making a rotation of coordinates. For a given $h>0$, by a rotation
of the coordinates we assume that $\inf\{|x|:\ x\in\Ga_{h, u}\}$ is
attained at $b^*e_n$, where $b^*\in (0, \ol b_h]$ and $e_k$ the unit
vector in the $x_k$-axis, $k=1, \cdots, n$. Then it suffices to
prove (2.34) for $\ol a_h, \ol b_h$ defined in this new coordinate
system. Since $\ol a_h\ge \ol b_h$, by the convexity of $\Ga_{h, u}$
we have
$$\Ga_{h,u}\cap\{|x|<\ol a_h\}\subset \{|x_n|\le 2\ol b_h\}. \tag 2.35$$

Let $\hat u$ be the restriction of $u$ on the 2-plane spanned by the
$x_1$ and $x_n$ axes. From the proof for the case $n=2$ in Lemma
2.1, we see that (2.34) holds if one can verify that $\Cal L_0 [\hat
u]\le C$ for some constant $C$ depending only on $n$, where $\Cal
L_0$ is the operator in (2.2).

For any given point $y=(y_1, 0, \cdots, 0, y_n)\in \Ga_{h, \hat u}$,
let $\kappa$ be the mean curvature of $\Ga_{h, u}$, and $\hat\kappa$
be the curvature of $\Ga_{h, \hat u}$ at $y$. Let $\gamma=(\gamma_1,
\gamma_2, \cdots, \gamma_n)$ be the unit normal of $\Ga_{h,u}$  at
$y$, and $\hat\gamma=(\hat \gamma_1, 0, \cdots, 0, \hat \gamma_n)$
be the unit normal of $\Ga_{h, \hat u}$  at $y$ in the 2-plane
spanned by the $x_1$ and $x_n$ axes. Suppose for a moment that
$$\gamma\cdot \hat\gamma
  =\gamma_1\hat\gamma_1+\gamma_n\hat\gamma_n\ge C_1\tag 2.36$$
for some positive constant $C_1$. Then by the convexity of $\Ga_{h,
u}$ we have $\hat \kappa\le C_2 \kappa$ and $\hat u_{\hat \gamma}\le
C_2u_\gamma$. By (2.2), we have $\kappa u_\gamma\le 1$. Hence $\Cal
L_0 [\hat u]=\hat\kappa \hat u_{\hat \gamma}\le C$ and so (2.34)
holds.

Let $P=\{x\in\R^n:\ \gamma\cdot (x-y)=0\}$ be the tangent plane of
$\Ga_{h, u}$ at the point $y$. Let $p_k=z_ke_k$, $k=1, \cdots, n$,
be the intersection of $P$ with the $x_k$-axis. Then by $\gamma\cdot
(x-y)=0$ at $x=p_k$, we have
$$\gamma_kz_k=\gamma_1y_1+\gamma_ny_n
 \ \ \ \forall\ \ k=2, \cdots, n-1$$
Hence if for all $k=2, \cdots, n-1$, $|z_k|\ge C\sqrt{y_1^2+y_n^2}$,
we have $|\gamma_k|\le \sqrt{\gamma_1^2+\gamma_n^2}$, which implies
$\sqrt{\gamma_1^2+\gamma_n^2}\ge C_1>0$ as $\gamma$ is a unit
vector. Observe that the vector $(\gamma_1, 0, \cdots, 0, \gamma_n)$
is parallel to the unit vector $\hat\gamma=(\hat \gamma_1, 0,
\cdots, 0, \hat \gamma_n)$. Hence we obtain (2.36).

To prove $|z_k|\ge C\sqrt{y_1^2+y_n^2}$, notice that $\gamma$ and
$\hat\gamma$ are invariant if we translate the level set $\Ga_{h,
u}$. Without loss of generality let us assume that $\gamma_n>0$. The
case $\gamma_n<0$ can be treated similarly. We translate $\Ga_{h,
u}$ in the $x_n$-direction by a distance $2\ol b_h$, so that
$\Ga_{h, u}\cap\{|x|<\ol a_h\}$ is contained in $\{x_n>0\}$. Since
$P$ is a tangent plane of $\Ga_{h, u}$ lying above the set $\Ga_{h,
u}$, we must have $|z_k|\ge \ol a_h$. On the other hand, $|y_1|\le
\ol a_h$ and $|y_n|\le 4\ol b_h$ (after the translation). Hence by
the assumption $\ol b_h\le \ol a_h$, we have $|z_k|\ge \frac
15\sqrt{y_1^2+y_n^2}$. \hfill $\square$

\proclaim{Lemma 2.7} Let $u$ be a complete convex solution of (1.2).
Suppose $u(0)=0$ and $u(\beta e_n)\ge 1$, $u(-\beta e_n)\ge 1$ for
some sufficiently small $\beta>0$, where $e_n=(0, \cdots, 0, 1)$.
Then $u$ is defined in a strip region.
\endproclaim

Note that the level set $\Ga_{h, u}=\{u=h\}$ may not be compact.
Note also that the strip region in Lemma 2.7 may not take the form
$\{x\in\R^n:\ -C_1\le x_n\le C_2\}$, except in some special cases
such as when $u$ is symmetric in $x_n$. But as in Lemmas 2.2 and
2.3, the axes of the minimum ellipsoid of $\Om_h\cap B_R(0)$ is a
small perturbation of axes of the coordinates.

To prove Lemma 2.7 we will prove that the graph of $g$, $\M_g=\{(x,
x_{n+1}): \ x_n=g(x', x_{n+1}), (x', x_{n+1})\in D\}$, is contained
in a strip $\{(x, x_{n+1})\in\R^{n+1}:\ 0\le x_n\le C\}$. By
convexity it suffices to prove $\ol b_h=g_h(0)$ is uniformly
bounded. The idea of our proof is very similar to the 2 dimensional
case given in \S2.1. In \S2.1 we divided the proof into three
lemmas. Here we present it in a single lemma.

\noo{\it Proof of Lemma 2.7.}\ \ Let $R=10^3$ and let $E$ be the
minimum ellipsoid of $\Om_1\cap B_R(0)$. By a rotation of
coordinates we assume the axial directions of $E$ coincide with
those of the coordinate system.

Denote $h_k=2^k$, $\ol a_k=\ol a_{h_k}$, $\ol b_k=\ol b_{h_k}$,
$g_k=g_{h_k}$, and $D_k=D_{h_k}$. As in \S2.1 we use induction
argument to prove
$$g_k (0)\le g_{k-1}(0)+ 2^{-k/4n}\
     \ \ \ \text{for all}\ \ k\ge 1. \tag 2.37$$
As shown in \S2.1, (2.37) implies that $u$ is defined in a strip
region.

The proof of (2.37) is similar to (2.3), we point out the difference
here. As in \S2.1, when $\beta$ is sufficiently small, by convexity
we have $\ol b_k\le h_k \ol b_0\le 2^k\beta\le 1$ when $k\le k_0$
and our induction argument starts at $k=k_0$.

Suppose by induction that (2.37) holds up to $k$. By the induction
assumption, $g_k(0)\le g_{k_0}(0)+\sum_{j=k_0}^k 2^{-j/4n}\le 2$. By
the concavity, $\ol b_{k+1}=g_{k+1}(0)\le 2g_k(0)\le 4$. Hence by
Lemma 2.6 we have
$$\ol a_{k+1}\ge C_0h_{k+1}.\tag 2.38$$

Next we prove (2.37) at $k+1$. Rotate the axes such that
$\p_ig_k(0)\le 0$ for all $i=1, \cdots, n-1$. By the concavity of
$g$ we have
$$g_k(0)=\sup\{g_k(x'): \ \ x_1>0, \cdots, x_{n-1}>0\}. \tag 2.39$$
Denote
$$L_k=\{x'\in\R^{n-1}:\ \ -\frac {C_0}{2n} h_k<x_i<\frac {C_0}{2n}h_k,
                               i=1, \cdots, n-1\}.$$
and $Q_k=L_k\times [h_k, h_{k+1}]\subset D$, where $C_0$ is the
constant in (2.38). Then similarly to (2.6),
$$\sup\{g(x',h):\ (x', h)\in Q_k\}
  \le 2\sup\{g(x', h_k):\ x'\in L_k\}
     \le 4 g(0, h_k)\le 8. $$
Observing that $2L_k\subset D_k$, by the convexity of $u$ we have
$L_k\subset D_{k-1}$. Hence by the concavity of $g$ and (2.39) we
have, for any $(x', h)\in Q_k$,
$$\align
|\p_h g(x', h)|
 & \le \frac {g(x', h)-g(x', h_{k-1})}{h-h_{k-1}} \tag 2.40 \\
 & \le \frac {g(x', h)}{h_k-h_{k-1}}
   \le \frac {2g(x', h_k)}{h_k-h_{k-1}}\\
 & \le \frac {2g(0, h_k)}{h_k-h_{k-1}}\le \frac 8{h_k}.\\
            \endalign $$
By (2.38) and (2.39), the concavity of $g$, and since $g\ge 0$,  we
also have,
$$|D_{x'} g(x', h)|\le C/h_k
 \ \ \ \forall\ (x', h)\in Q_k.\tag 2.41$$

From the above gradient estimates and the concavity of $g$, the
average in $Q_k$ of the second order derivatives $|\p^2 g|\le
Ch_k^{-2}$. But we have not proved this estimate pointwise. We need
to treat the set of points where $|\p^2 g|$ is relatively large.
Denote
$$\align
\chi & = \{(x', h)\in Q_k:
   |{\Sigma}_{i=1}^{n-1} \p_i^2 g_h(x')| \ge h_k^{-5/4}\} , \tag 2.42\\
\chi^+ & = \{(x', h)\in Q_k:
   |{\Sigma}_{i=1}^{n-1} \p_i^2 g^+_h(x')| \ge h_k^{-5/4}\} .\\
   \endalign    $$
Obviously $\chi^+\subset\chi$. By the gradient estimates, we have
$$|\chi|_{\Cal H^n} h_k^{-5/4}
   \le \big{|}\int_{Q_k}\Delta_{x'} g\big{|}
   \le \int_{\p L_k\times [h_k, h_{k+1}] } |D_{x'} g|
   \le Ch_k^{n-2}.$$
In the above formula $g$ is a function of $(x', h)$.  Hence we
obtain
$$|\chi^+|_{\Cal H^n}\le |\chi|_{\Cal H^n}\le C h_k^{n-3/4} . \tag 2.43$$

From (2.43) and by the Fubini Theorem, there is a set $\wtt L\subset
L_k$ with measure $|\wtt L|_{\Cal H^{n-1}} < h_k^{n-3/2}$ such that
for any $y'\in L_k-\wtt L$,
$$|\chi^+_{y'}|_{\Cal H^1}\le C h_k^{3/4}, \tag 2.44$$
where $\chi^+_{y'}=\chi^+\cap\{x'=y'\}$.

For any given $y'\in L_k-\wtt L$, we want to prove
$$g^+_{k+1}(y')-g^+_k(y')\le Ch_k^{-1/4}.\tag 2.45$$
Similarly we can estimate $|g^-_{k+1}(y')-g^-_k(y')|$. Hence if
(2.45) is proved, we have
$$g_{k+1}(y')-g_k(y')\le Ch_k^{-1/4},$$
which corresponds to (2.15).  As the argument after (2.15), we can
choose a point $y_1\in L_k-\wtt L$ with $|y_1|\le Ch_k^{1-1/2(n-1)}$
such that $g_{k+1}(y_1)\le g_k(0)+Ch_k^{-1/4}$. But now
$$\frac{g_{k+1}(0)}{g_{k+1}(y_1)}
 \le \frac {\ol a_{k+1}}{\ol a_{k+1}-|y_1|}\le 1+Ch_k^{-1/2(n-1)}. $$
Therefore we obtain (2.37).

To prove (2.45), we have
$$\align
g^+_{k+1}(y')-g^+_k(y')
 &=\int_I\p_h g^+(y', h) dh\tag 2.46\\
 &=\left(\int_{I_1}+\int_{I_2}+\int_{I_3}\right) \p_h g^+(y', h) dh,\\
 \endalign $$
where as in (2.24), $I=(h_k, h_{k+1})$, $I_1= \chi^+\cap\{x'=y'\}$,
$I_2 = \{h\in I-I_1:\ \frac{\sigma
u_{\gamma\gamma}(q)}{\sigma+u_\gamma^2(q)} \le \frac 12\}$, and $
I_3 =I-(I_1\cup I_2)$, where $q=\Cal T(p)$ with $p=(y', h)$, and
$\Cal T: p\to q$ is the mapping introduced before (2.21).

For the first integral in (2.46), by (2.40) and (2.44) we have
$$\int_{I_1} \p_hg^+(y', h)dh \le Ch_k^{-1/4}$$
Note that in $I_2$, similarly to (2.23) we have
$$|\p_h g^+(y_1, h)|\le C|\p_{x_1x_1}g^+|\le Ch_k^{-5/4}
 \ \ \forall\ (y_1, h)\not\in \chi^+. $$
Hence we have the estimate for the second integral in (2.46),
$$\int_{I_2} \p_h g^+(y', h)dh\le  Ch_k^{-1/4}. $$
For the third one, the argument between (2.25) and (2.29) applies
and we also have the estimate $\text{osc}_{I_3}g^+\le C/h_k$. Hence
(2.45) holds. \hfill$\square$

\vskip5pt

The next lemma corresponds to Lemma 2.4 in \S2.1.

\proclaim{Lemma 2.8} Let $u$ be a complete convex solution of (1.2).
Suppose $u(0)=0$ and the infimum $\inf\{|x|:\ x\in \Ga_1\}$ is
attained at $x_0=(0, \cdots, 0, -\delta)\in \Ga_1$ for some
$\delta>0$ sufficiently small. As above let $D_1$ be the projection
of $\Ga_1$ on the plane $\R^n\cap \{x_n=0\}$. Then $D_1\supset
\{x'\in\R^{n-1}:\ |x'|<R\}$ with
$$R\ge \frac 1{C_n} (-\log\delta-C)^{1/2},\tag 2.47$$
where $C_n$ is a constant depending only on $n$, and $C>0$ is a
constant independent of $\delta$.
\endproclaim

\noo{\it Proof.}\ \ Estimate (2.47) is equivalent to $a_1(p)\ge
\frac 1{C_n} (-\log\delta-C)^{1/2}$ for any $p\in S^{n-2}$. Suppose
$\inf a_1(p)$ is attained at $p=(1, 0, \cdots, 0)$. By restricting
$u$ to the 2-plane $\{x_2=\cdots=x_{n-1}=0\}$, we reduce the proof
to the 2 dimensional case in Lemma 2.4, as we have shown, in the
proof of Lemma 2.6, that $\Cal L_0 [u]\le C_n$. \hfill $\square$

With Lemmas 2.7 and 2.8, we extend Theorem 2.1 to high dimensions.

\proclaim{Theorem 2.2} Let $u$ be an entire convex solution of (1.2)
in $\R^n$. Then there exists a positive constant $C$ such that for
any $x\in\R^n$,
$$u(x)\le C (1+|x|^2),  \tag 2.48$$
where $C$ depends on $n$ and the upper bound of $u(0)$ and
$|Du(0)|$.
\endproclaim

\noo{\it Proof.} The proof is very similar to that of Theorem 2.1.
Let $\delta, \delta^*$ and $u_h$ be as in the proof of Theorem 2.1.
Instead of an ellipse, here we use the ellipsoid
$$E=\{x\in\R^2:\ {\sum}_{i=1}^{n-1}\frac {x_i^2}{(R/2n)^2}
   +\frac {|x_2-(\delta^*-\delta)/2|^2}
   {[(\delta^*+\delta)/8]^2}<1\}.$$

When $\sigma=0$, the level set $\Gamma_{-t, u_h}$ is a solution to
the mean curvature flow. Let $E_{-t}$ be the solution to the mean
curvature flow with initial condition $E_{-1}=E$, so that
$E_{-t}\subset \Om_{-t, u_h}$ for all $t>-1$. Suppose it takes time
$T$ for $E_{-t}$ to shrink to a point. Then we have $\inf_{B_1(0)}
u_h\le 1-T$. Observe that for any fixed $\delta^*$, $E_{-t}$
converges to a pair of parallel planes, and so $T\to\infty$ as
$R\to\infty$ (or $\delta\to 0$). Hence when $\delta$ is small, we
reach a contradiction with the assertion that $\inf_{B_1(0)} u_h$ is
uniformly bounded from below. The case $\sigma>0$ can be proved in
the same way as in Theorem 2.1. \hfill$\square$

\vskip5pt

From Theorem 2.2,  we have accordingly

\proclaim{Corollary 2.2} Let $u$ be a complete convex solution of
(1.2). Then $u$ is either an entire solution, or is defined in a
strip region. There is no complete convex solution of (1.2) defined
in a half space.
\endproclaim

Note that estimate (2.48) also implies the follows compactness
result. This compactness result is not just for the set of blow-up
solutions to mean convex flow but for all entire convex solutions of
(1.2). We don't know whether an entire convex solution to (1.2) must
be a blow-up solution to mean convex flow.

\proclaim{Corollary 2.3} For any constant $C>0$, the set of all
entire convex solutions $u$ to (1.2) satisfying $u(0)=0$ and
$|Du(0)|\le C$ is compact.
\endproclaim

\vskip10pt

\noo{\bf 2.3. Proof of Theorem 1.3}. First we prove a lemma.

\proclaim{Lemma 2.9} Let $u$ be an entire convex solution of (1.2).
Suppose $u\ge 0$ and $u(0)=0$. Then the convex set $\{u=0\}$ is
either a single point or it is a linear subspace of $\R^n$.
\endproclaim

\noo{\it Proof.}\ If $\sigma>0$, $u$ is analytic. As the set
$\{u=0\}$ is convex, it must be a single point or a linear
subspace of $\R^n$. In the following we consider the case
$\sigma=0$.

If the set $\{u=0\}$ is bounded, then $\Ga_{h, u}$ is a closed,
bounded convex hypersurface. As $\Ga_{h, u}$ evolves by mean
curvature (with time $t=-h$). From [8, 12] it follows that
$\{u=0\}$ is a single point.

If the set $\{u=0\}$ contains a straight line, say the line
$\ell=(t, 0,\cdots, 0)$ ($t\in\R$), then by convexity $u$ is
independent of $x_1$. Hence to prove Lemma 2.9, we need only to rule
out the possibility that $\{u=0\}$ contains a ray but no straight
line lies in it.

Suppose the ray $r=(t, 0, \cdots, 0)$ ($t>0$) is contained in
$\{u=0\}$. We may also suppose that $\{u=0\}$ contains no straight
lines and the asymptotical cone of $\{u=0\}$ is contained in
$\{x_1>0\}$. Then $u$ is decreasing in $x_1$. Denote $u_m (x_1,
x_2, \cdots, x_n)=u(x_1+m, x_2, \cdots, x_n)$, where $m>0$ is a
constant. Then $u_m$ is nonnegative and decreasing in $m$. By
choosing a subsequence we suppose $u_m\to \hat u$ as $m\to\infty$.
Then the straight line $\ell=(t, 0,\cdots, 0)$ ($t\in\R$) is
contained in the graph of $\hat u$. By convexity, $\hat u$ is
independent of $x_1$. Since $\Cal L_0[\hat u]=1$, $\hat u$ does
not vanish completely, and so we must have $n\ge 3$. Moreover, we
have $\hat u<u$ except on the set $\{u=\hat u=0\}$.

Since $u$ and $\hat u$ are both solutions to $\Cal L_0[u]=1$, the
level sets $\{u=-t\}$ and $\{\hat u=-t\}$ evolve by mean curvature
(with time $t$). Denote $\M_t=\{u=-t\}\cap\{x_1=0\}$ and
$\hat\M_t=\{\hat u=-t\} \cap \{x_1=0\}$. Then $\hat \M_t$ evolves
by mean curvature as $\hat u$ is independent of $x_1$. We assert
that $\M_t$ evolves at a velocity greater than its mean curvature.
Indeed, for any given point $p\in\M_t$, we assume the hypersurface
$\{u=-t\}$ is locally given by $x_n=\psi(x_1, \cdots, x_{n-1})$,
and locally $\M_t$ is given by $x_n=\psi(0, x_2, \cdots,
x_{n-1})$.  By choosing the coordinate system properly we also
assume that $\p_{x_i}\psi=0$ for $i=2, \cdots, n-1$ at $p$. Then
$\M_t$ evolves at the velocity $\sqrt{1+|D\psi|^2}\, \div \frac
{D\psi}{\sqrt{1+|D\psi|^2}}$, by convexity which is greater than
$\sum_{i=2}^{n-1}\p_{x_i}^2\psi$, the mean curvature of $\M_t$ at
$p$.

On the other hand, since $\hat u<u$, $\M_t$ is strictly contained in
the interior of $\hat\M_t$ for any $t<0$. Moreover $\M_t$ is a
bounded, closed convex hypersurface, as the asymptotical cone of
$\{u=0\}$ is contained in $\{x_1>0\}$. By the comparison principle,
$\M_t$ is strictly contained in the interior of $\hat\M_t$ for all
$t\le 0$. We reach a contradiction as $\hat u=u=0$ at the origin.
\hfill$\square$

\vskip5pt

Therefore by [12], Lemma 2.9 implies that the singularity set of a
mean curvature flow of convex, noncompact hypersurfaces in
$\R^{n+1}$ must be a subspace $\R^{n-k}$ for some $1\le k\le n$. But
then by convexity, $u$ is a function of $k$ variables.

\vskip10pt

\noo{\it Proof of Theorem 1.3}.\ \
\newline
{\it Step 1}. First we prove that there is a subsequence of $u_h$,
where $u_h(x)=h^{-1}u(h^{1/2}x)$, which converges to $\eta_k$ for
some $2\le k\le n$, where $\eta_k$ is the function given in (1.4).

By adding a constant we may suppose $u(0)=0$. Let $T=\{x_{n+1}=
a(x)\}$ be the tangent plane of $u$ at the origin. By Theorems 2.1
and 2.2 and the convexity of $u$ we have
$$a(x)\le u(x)\le C(1+|x|^2).$$
Hence
$$\frac {1}{\sqrt h} a(x) \le u_h(x)
      \le C(\frac {1}{\sqrt h}+|x|^2).$$
By convexity it follows that $Du_{h}$ is locally uniformly
bounded. Hence $u_h$ sub-converges to a convex function $u_0$
which satisfies $u_0(0)=0$,
$$0\le u_0(x)\le C|x|^2.\tag 2.49$$
By Lemma 2.5, $u_0$ is an entire convex solution of $\Cal
L_0[u]=1$.

Case 1: the set $\{x\in\R^n:\ u_0(x)=0\}$ is bounded. Then by
convexity the level set $\Ga_{1, u_0}=\{x\in\R^n:\ u_0(x)=1\}$ is
a bounded convex hypersurface. Since the level set $\{u_0=-t\}$,
with time $t\in (-\infty, 0)$, evolves by mean curvature, by the
asymptotic estimates in [8, 12],
$$u_0(x)=\frac {1}{2(n-1)}|x|^2+\phi(x)\tag 2.50$$
where $\phi(x)=o(|x|^2)$ for $x\ne 0$ near the origin. Hence for
any $\eps>0$, there is a sufficiently small $h'>0$, such that
$$B_{(1-\eps)r}(0)\subset \Om_{h',u_0}\subset B_{(1+\eps) r}(0), $$
where $r=\sqrt{2(n-1)h'}$. Hence there is a sequence $h_m\to
\infty$ such that
$$B_{(1-\frac 1m)r_m}(0)\subset \Om_{h_m,u}
     \subset B_{(1+\frac 1m)r_m}(0), \tag 2.51$$
where $r_m=\sqrt{2(n-1)h_m}$. Let  $u_{h_m}(x)=\frac 1{h_m} u(\sqrt
{h_m}\, x)$. Then $u_{h_m}$ sub-converges to $\hat u_0$ which
satisfies $\Cal L_0[\hat u_0]=1$. From (2.51), the level set
$\Ga_{1, \hat u_0}$ is a sphere. Hence $\hat u_0(x)=\frac
{1}{2(n-1)}|x|^2$.

Case 2: the set $\{u_0=0\}$ is unbounded. Then by Lemma 2.9, the set
$\{u_0=0\}$ is a linear sub-space of $\R^n$. Suppose $\{u_0=0\}=
\{x\in\R^n:\ x_{k+1}=\cdots=x_n=0\}$. We must have $k\ge 2$, as the
level set $\{u_0=-t\}$ evolves by its mean curvature. It follows
that $u_0$ is a convex function depending only on $\hat x=(x_1,
\cdots, x_k)$. Similarly as above we have
$$u_0(x)=\frac {1}{2(k-1)}|\hat x|^2 +o(|\hat x|^2) \tag 2.52$$
near the origin. Hence for any $\eps>0$,
$$\hat B_{(1-\eps)r}(0)\subset \hat\Om_{h',u_0}
     \subset \hat B_{(1+\eps)r}(0) $$
provided $h'$ is sufficiently small, where $r=\sqrt{2(k-1)h'}$,
$\hat B_r(0)=B_r(0)\cap\{\wtt x=0\}$ and $\hat \Om_{h', u}=\Om_{h',
u}\cap \{\wtt x=0\}$, $\wtt x=(x_{k+1}, \cdots, x_n)$. It follows
that for any $R>0$,
$$\align
\{x\in \R^n:\ |\hat x|<(1-\eps)r\}\cap\{|\wtt x|<R\}
  & \subset \Om_{h', u_{h_m}}\cap\{|\wtt x|<R\}\\
 & \subset \{x\in \R^n:\ |\hat x|<(1+\eps)r\}\cap\{|\wtt x|<R\}\\
 \endalign $$
if $h_m$ is sufficiently large. Hence there exist $\tau_m\to\infty$
and (a different sequence) $h_m\to\infty$ such that
$$\align
\{x\in \R^n:\ |\hat x|<(1&-\frac 1m)r_m\}\cap\{|\wtt x|<\tau_m r_m\}
                                           \tag 2.53\\
 &\subset \Om_{h_m, u} \cap \{|\wtt x|<\tau_m r_m\} \\
 &\subset \{x\in \R^n:\ |\hat x|<(1+\frac 1m)r_m  \}
                     \cap \{|\wtt x|<\tau_m r_m\},\\
 \endalign $$
where $r_m=\sqrt{2(k-1)h_m}$. Hence $u_{h_m}\to \frac
{1}{2(k-1)}|\hat x|^2$.

\vskip10pt

\noo {\it Step 2}.\ Now we prove that $u_h$ itself, after a rotation
of axes, converges to the function $\eta_k$.

In Step 1 we proved that $u_{h_m}$ converges to $\eta_k$ for some
$2\le k\le n$. Let us choose the sequence $\{h_m\}$ properly such
that $k$ is the largest such integer, namely if $u_{h_m'}$ converges
to $\eta_{k'}$, then $k'\le k$. From the above proof we can also
choose $h_m$ such that (2.51) or (2.53) holds.

Case 1: $k=n$. We prove that for any constant $\eps>0$,
$$B_{(1-\eps)r}(0)\subset \Om_{h,u}
     \subset B_{(1+\eps)r}(0) \tag 2.54$$
if $h>0$ is sufficiently large, where $r=\sqrt{2(n-1)h}$. Suppose
(2.54) is not true. Let $h_m\to\infty$ be a sequence such that
(2.51) holds. Let
$$\hat h_m=\inf\{h'\le h_m:
        \ (2.54)\ \text{holds for any}\ h\in (h',h_m)\}. \tag 2.55$$
Since $u_{h_m}\to \frac {1}{2(n-1)}|x|^2$, we have $h_m/\hat h_m
\to\infty$ as $m\to\infty$. Let $\alpha>1$ be a fixed constant which
will be determined below. Then the sequence $u_{\alpha \hat h_m}$
sub-converges to a convex function $u_0$ satisfying $u_0(0)=0$,
$u_0\ge 0$ and $\Cal L_0[u]=1$.  By our choice of $\hat h_m$, the
level set $\Om_{1, u_0}$ satisfies
$$B_{(1-\eps)r}(0)\subset \Om_{1,u_0}
     \subset B_{(1+\eps)r}(0)$$
with $r=\sqrt{2(n-1)}$. Hence from [8, 12], the level set $\Om_{h,
u_0}$ satisfies
$$B_{(1-\delta)r}(0)\subset \Om_{h,u_0}
     \subset B_{(1+\delta)r}(0)$$
with $\delta\to 0$ as $h\to 0$, where $r=\sqrt{2(n-1)h}$. Hence we
have
$$B_{(1-2\delta)r}(0)\subset \Om_{h,u_{\alpha \hat h_m}}
     \subset B_{(1+2\delta)r}(0)$$
if $m$ is sufficiently large. Choose $h$ sufficiently small such
that $\delta\le\frac 13 \eps$ and let $\alpha=h^{-1}$. Then scaling
back we find that $\Om_{\hat h_m, u}$ satisfies
$$B_{(1-2\delta)r_m}(0)\subset \Om_{\hat h_m,u}
     \subset B_{(1+2\delta)r_m}(0)$$
with $r=\sqrt{2(n-1)\hat h_m}$. When $\delta<\frac 12\eps$, this is
in contradiction with our choice of $\hat h_m$.

Case 2: $k<n$. For any given small $\eps>0$, by (2.53), $\Ga_{h_m,
u}$ is $\eps$-close to the cylinder $S^{k-1}\times R^{n-k}$ if $m$
is sufficiently large, namely
$$\align
\{x\in \R^n:\ |\hat x|<&(1-\eps)r\}\cap\{|\wtt x|<\eps^{-1} r\} \tag 2.56\\
 &\subset \Om_{h_m, u}\cap\{|\wtt x|<\eps^{-1} r\}\\
 &\subset \{x\in \R^n:\ |\hat x|<(1+\eps)r\}\cap\{|\wtt x|<\eps^{-1} r\},\\
 \endalign $$
where $r=\sqrt{2(k-1)h_m}$. Let $\hat h_m<h_m$ be the least number
such that $\Ga_{h, u}$ is $\eps$-close to the cylinder
$S^{k-1}\times R^{n-k}$ (the axes of the cylinder may vary as $h$
varies) for any $h\in [\hat h_m, h_m]$. Then by our assumption that
$k$ is the largest possible integer, we have, due to (2.52), that
$u_{\alpha \hat h_m}= \frac {1}{2(k-1)}|\hat x|^2+o(|\hat x|^2)$ for
any given $\alpha>1$. Here we regard $u_{\alpha \hat h_m}$ as a
function of $\hat x=(x_1, \cdots, x_k)$ by letting
$x_{k+1}=\cdots=x_n=0$. Similar to Case 1, we can choose $\alpha>1$
such that $\Ga_{\hat h_m, u}$ is $\frac 12\eps$-close to the
cylinder $S^{k-1}\times R^{n-k}$, which is in contradiction with our
choice of $\hat h_m$. Hence Theorem 1.3 is proved. $\square$

\vskip5pt

Note that Case 2 in Step 2 follows readily from Step 1 and the fact
that $k$ is an integer. In Step 1 it is shown that $u_h$ converges
along a subsequence to the function $\eta_k$ but $k$ is an integer
so it must be the same integer for all subsequences.

\vskip5pt

\noo{\bf Remark 2.2.} \ Theorem 1.3 asserts that $u_h$, which is the
blow-down of $u$ with respect to the origin in space-time,
sub-converges to a self-similar solution. We point out that for
ancient convex solutions $w$ to the level set flow (1.2) (with
$\sigma=0$), under some very mild conditions the corresponding level
set $\{w=h\}$, after proper translation, sub-converges as
$h\to\infty$ to a translating solution. In particular, if $w$ is a
complete convex solution of (1.2) defined in a strip region in
$\R^2$, then after proper translation, the level set must converge
along a subsequence to the grim reaper.

 \vskip20pt

\centerline{\bf 3. Proof of Theorem 1.1}

\vskip10pt

In this section we prove Theorem 1.1. Throughout this section we
suppose the dimension $n=2$.

Let $u$ be an entire convex solution of (1.2). By Theorem 1.3 we
have
$$u(x)=\frac 12 |x|^2 +\phi(x)\tag 3.1$$
with $|\phi(x)|=o(|x|^2)$ as $|x|\to\infty$. To prove Theorem 1.1,
we first consider the case $\sigma=0$.

\proclaim{Theorem 3.1} Let $u$ be an entire convex solution of
(1.2) with $\sigma=0$.  Then $u(x)=\frac 12 |x|^2$ in a proper
coordinate system.
\endproclaim

\noo{\it Proof}.\ \ \  By a translation of the graph of $u$, we may
suppose $u\ge 0$, $u(0)=0$, and (3.1) holds. For any constant $h>1$,
denote $u_h(y)=u(h^{1/2}y)/h$. Then $u_h$ is also an entire convex
solution of (1.2) and by (3.1), the sub-level set $\Om_{1/2, u_h}$
satisfies
$$B_{1-\eps}(0)\subset \Om_{1/2,u_h}
     \subset B_{1+\eps}(0) \tag 3.2$$
with $\eps\to 0$ as $h\to\infty$. By Gage-Hamilton [8], we have
$$u_h(y)=\frac {1}{2}|y|^2 +\phi(y)\tag 3.3$$
with
$$|\phi (y)|\le C|y|^{2+\alpha}$$
for some $\alpha\in (0, 1)$, and $C$ is a constant independent of
$h$. Rescaling back to the $x$-coordinate we obtain
$$u(x)=\frac {1}{2}|x|^2 + h^2 \phi(x/h) ,$$
where for any fixed $x$, $h^2 \phi(x/h)\to 0$ as $h\to\infty$. Hence
$u(x)\equiv \frac {1}{2}|x|^2$. \hfill$\square$

\noo{\bf Remark 3.1}.\ By the asymptotic estimates in [12], Theorem
3.1 also holds in high dimensions if the solution $u$ satisfies
$$C_1|x|^2\le u(x)\le C_2 |x|^2 .\tag 3.4$$
Indeed, if $u$ satisfies (3.4), we have $u(x)=\frac {1}{2(n-1)}
|x|^2 +o(|x|^2)$ by Theorem 1.3. Next we consider the case
$\sigma=1$ of Theorem 1.1.

\proclaim{Theorem 3.2} Let $u$ be an entire convex solution of the
mean curvature equation (1.1).  Then $u$ is rotationally symmetric
in an appropriate coordinate system.
\endproclaim

To prove Theorem 3.2 we need a few lemmas.

\proclaim{Lemma 3.1} Let $\Om$ be a bounded convex domain in
$\R^2$. Let $u_0$ and $u_\sigma$ be respectively solutions of
$\Cal L_0[u]=1$ and $\Cal L_\sigma[u]=1$ in $\Om$, vanishing on
$\pom$, where $\sigma\in (0, 1]$. Suppose $u_\sigma$ is convex.
Then for any constant $a>0$, there exists a constant $C>0$,
depending on $a$ and the upper and lower bounds of $|Du_\sigma|$
on the set $\{x\in\Om:\ \inf u_\sigma+a\le u_\sigma(x)<0\}$, such
that for any $0>h>a+\inf_\Om u_\sigma$,
$$0\le |\Om_{h, u_0}|-|\Om_{h, u_\sigma}|\le C\sigma.\tag 3.5$$
\endproclaim

Note that the constant $C$ in (3.5) is large when the lower bound of
$|Du_\sigma|$ is small. We don't impose condition on $\Om$ but it is
convex and its shape is controlled by the lower bound of
$|Du_\sigma|$.

For the proof of Theorem 1.1, the solution $u$ satisfies (3.1) and
$\Om$ is a small perturbation of the unit disc. In this case (3.5)
can be proved easily. In fact, the difference $|\Om_{h,
u_0}|-|\Om_{h, u_\sigma}|$ is controlled by $\frac {\sigma
u_{\gamma\gamma}}{\sigma+u_\gamma^2}=O(\sigma)$, see (3.6) below.

\vskip5pt

\noo{\it Proof}.\ \  Denote $u=u_\sigma$ and suppose without loss
of generality that $u(0)=\inf_\Om u=-1$. By convexity we have
$\L_0[u]\le 1$. By the comparison principle we have $u\ge u_0$.
Hence $\Om_{h, u}\subset\Om_{h, u_0}$ and $|\Om_{h, u_0}|\ge
|\Om_{h, u}|$.

Write the equation $\Cal L_\sigma[u]=1$ in the form
$$ \kappa u_\gamma
  =1-\frac {\sigma u_{\gamma\gamma}}{\sigma+u_\gamma^2},\tag 3.6$$
where $\kappa$ is the curvature of the level set $\Gamma_{h, u}$,
and $\gamma$ is the unit outward normal to $\Om_{h, u}$. (3.6)
implies that the level set $\Gamma_{h, u}$ is moving  with the
velocity (regard $t=-h$ as the time)
$$v=u_\gamma^{-1}=\frac {\kappa}
    {1-\frac {\sigma u_{\gamma\gamma}}{\sigma+u_\gamma^2}}.$$

Let $w=w(\cdot, h)\in C (S^1)$ denote the supporting function of
$\Gamma_{h, u}$. That is,
$$w(p)=w(p, h)
   =\sup\{\lan p, x\ran:\ x\in \Gamma_{h, u}\}\ \ \ p\in S^1. $$
The supremum is attained at the point $x$ at which the unit outer
normal $\gamma(x)=p$, and the curvature $\kappa$ at $x$ is given
by
$$\kappa(x)=\frac {1}{(w''+w)(p)},$$
where $S^1$ is parametrized by $p=(\cos\th,\sin\th)$ and $w'=\frac
{d}{d\th}w$. The area of the domain $\Om_{h, u}$ is given by
$$|\Om_{h, u}|=\frac 12 \int_{S^1} w(w''+w).$$
Observing that $\p_h w=v=u_\gamma^{-1}$, we have
$$\align
\frac {d}{dh}|\Om_{h, u}|
 &= \frac {d}{dh}\int_{S^1}\frac 12 w(w''+w)\tag 3.7\\
 &= \int_{S^1} \p_h w(w''+w)\\
 &= \int_{S^1} \frac {\kappa}
  {1-\frac {\sigma u_{\gamma\gamma}}{\sigma+u_\gamma^2}}(w''+w).\\
 \endalign $$

For any $h\in (a+\inf_\Om u, 0)$, denote $D =S^1 \times (h, 0)$.
Let $G$ denote the diffeomorphism from $\M_{u, h}=:
\M_u\cap\{h<u<0\}$ to $D$, where $\M_u$ is the graph of $u$, such
that for any point $(x, t)\in\M_{u, h}$,  $G(x, t)=(G_t(x), t)\in
D$, where $G_t$ is the Gauss mapping from the level set
$\Gamma_{t,u}$ to $S^1$.

We divide $D$ into two parts, $D = D_1\cup D_2$, such that
$$D_1=\{1-\frac {\sigma u_{\gamma\gamma}}{\sigma + u_\gamma^2}
   \ge \frac 12\}$$
and $D_2=D-D_1$. Observing that $\kappa (w''+w)=1$, from (3.7) we
have
$$|\Om_{0, u}|-|\Om_{h, u}|
   =\int_{D_1} \frac {1}
     {1-\frac {\sigma u_{\gamma\gamma}}{\sigma+u_\gamma^2}}
    +\int_{D_2}\frac {\kappa}
  {1-\frac {\sigma u_{\gamma\gamma}} {\sigma+u_\gamma^2}}(w''+w).$$
On $D_1$ we have
$$(1-\frac {\sigma u_{\gamma\gamma}}{\sigma+u_\gamma^2})^{-1}
 \le 1+\frac {2\sigma u_{\gamma\gamma}}{\sigma+u_\gamma^2}
 \le 1+C_1\sigma u_{\gamma\gamma},$$
on $D_2$ we have
$$  \frac {\kappa}
  {1-\frac {\sigma u_{\gamma\gamma}}{\sigma+u_\gamma^2}}
       = u_\gamma^{-1}\le C_2^{-1}, \tag 3.8$$
where both constants $C_1$ and $C_2$ depend on the lower bound of
$|Du|$ on the set $\{u>a+\inf u\}$.  Hence
$$\align
|\Om_{0, u}|-|\Om_{h, u}|
 &\le  \int_{h}^0\int_{S^1} (1+C\sigma u_{\gamma\gamma})
    +C\int_{D_2}(w''+w)  \tag 3.9\\
  & = 2\pi |h|+C\sigma \int_{h}^0\int_{S^1} u_{\gamma\gamma}
    +C|G^{-1}(D_2)|\\
  & \le 2\pi|h|+C\sigma+C|G^{-1}(D_2)|\\
\endalign$$

To estimate $|G^{-1}(D_2)|$ we suppose $\inf u$ is attained at the
origin. For any unit vector $\tau$ in the plane $\{x_3=0\}$ starting
at the origin, let $P_\tau$ be the plane in $\R^3$ containing $\tau$
and the $x_3$-axis, and let $E_\tau$ denote the intersection of
$P_\tau$ with $G^{-1}(D_2)$. On $E_\tau$ we have, by our definition
of $D_2$, $u_{\gamma\gamma}\ge \frac {C}{\sigma}$. Noting that by
equation, $u_{\xi\xi}\le C$ in $G^{-1}(D_2)$ for any unit vector
tangential to $\Ga_{h, u}$ and that the inner product $\lan \gamma,
\tau\ran\ge C'$ for some constants $C, C'>0$ depending on the upper
and lower bounds of $|Du|$ on the set $\{x\in\Om:\ \inf u +a\le
u(x)<0\}$ (which also determine the geometric shape of $\Om$), we
have $u_{\tau\tau}\ge \frac {C}{\sigma}$ for a different $C$ (for
small $\sigma>0$). It follows that the one dimensional Lebesgue
measure $|E_\tau|_{\Cal H^1}\le C\sigma$ for some $C$ depending on
the upper bound of $|Du|$. Hence the two dimensional Lebesgue
measure $|G^{-1}(D_2)|_{\Cal H^2}\le C\sigma$.

Observing that the level set $\Gamma_{h, u_0}$ is moving by its
curvature (with time $t=-h$), we have
$$|\Om_{0, u_0}|-|\Om_{h, u_0}| =2\pi |h|.$$
Hence by (3.9),
$$ |\Om_{h, u_0}|-|\Om_{h, u}|
  \le C\sigma +C|G^{-1}(D_2)|
 \le C\sigma .$$
We obtain (3.5). \hfill$\square$

\proclaim{Lemma 3.2}  Let $\{\ell_t\}$ be a convex solution to the
curve shortening flow. Suppose $\ell_0$ is in the
$\delta_0$-neighborhood of a unit circle $S^1$ and $\{\ell_t\}$
shrinks to a point (the origin) at $t=\frac 12$. Let $\hat\ell_t=
\frac {1}{\sqrt{1-2t}}\ell_t$ be the normalization of $\ell_t$. Then
$\hat\ell_t$ is in the $\delta_t$-neighborhood of the unit circle
centered at the origin,
$$\hat \ell_t\subset N_{\delta_t}(S^1) , \tag 3.10$$
with
$$\delta_t\le C\delta_0 (\frac 12-t)^\alpha, $$
where $\alpha\in (0, 1)$ is a positive constant.
\endproclaim

\noo{\it Proof}.\ \ First observe, by the comparison principle,
that when $t\le \frac 14$, $\ell_t$ is pinched between two
concentrated circles with Hausdorff distance $C\delta_0$. By the
Schauder estimate, for $t\in (\frac 18, \frac 14)$ the $C^k$ norm
of $\hat \ell_t$ is in the $C\delta_0$-neighborhood of the unit
circle, that is
$$\|\hat\ell_t-S^1\|_{C^k}\le C\delta_0.\tag 3.11$$
With estimate (3.11) we obtain (3.10) from [8], \S 5.7.10-\S 5.7.15.
\hfill$\square$

\noo {\bf Remark 3.2 (i)}.  By the Schauder estimate one can
simplify some estimates in [8], \S 5.1-5.6.  In [8], \S 5.7.10-\S
5.7.15, it was proved that for any $\alpha>0$ small, there exists
$\delta_0>0$ such that if (3.11) holds at $t=0$, then
$$\int_{S^1} [\kappa'(\tau)]^2
     \le e^{-\alpha \tau }\int_{S^1} [\kappa' (0)]^2,\tag 3.12$$
where $\tau=\frac 12\log (\frac 12-t)$, and $\kappa'$ denotes the
derivative of the curvature $\kappa$ with respect to the
are-length parameter. Similar inequalities for high order
derivatives of $\kappa$ were also proved there.

(ii) Let $u$ be a convex solution of $L_0[u]=0$ which attains its
minimum $0$ at $y_1$ (namely $u(y_1)=\inf u=0$). Suppose the level
set $\Ga_{1/2}\subset N_{\delta_0}(S^1)$ for some small
$\delta_0>0$. Then $|y_1|<C\delta_0$ for some $C>0$ independent of
$\delta_0$. Therefore by a translation we may assume that $u$
attains its minimum at $0$ and $\Ga_{1/2, u}\subset
N_{C^*\delta_0}(S^1)$ for a different constant $C^*$.

To prove $|y_1|<C\delta_0$, let $\hat u=\frac 12 |x|^2$ be the
rotationally symmetric solution to $L_0[u]=0$. As in the proof of
Lemma 3.2, let $w(p, h)$ and $\hat w(p, h)$ be respectively the
support functions of $\Ga_{h, u}$ and $\Ga_{h,\hat u}$, where
$p=(\cos\th, \sin\th)$. Denote $t=-h$ (regard $t\in (-\frac 12, 0)$
as the time). Then $w_t(w''+w)=-1$, $\hat w_t(\hat w''+\hat w)=-1$.
Denote $\phi=w-\hat w$. Direct computation shows that
$$(w''+w+\hat w''+\hat w)\phi_t=-(w_t+\hat w_t)(\phi''+\phi).$$
Hence $\phi$ satisfies the equation
$$\phi_t =(w_t\hat w_t)(\phi''+\phi). $$
We have $\hat w_t=-\frac{1}{\sqrt {|t|}}$ and by estimate (3.12)
(for higher order derivatives), $w_t=-\frac {1+o(1)}{\sqrt {|t|}}$,
as $t\to 0$. We obtain $w_t\hat w_t=\frac {1+o(1)}{|t|}$. The
estimate (3.12) also implies that the curvature of $\Ga_{u, h}$ is
equal to that of $\Ga_{\hat u, h}$ up to a lower order perturbation,
namely $\phi''+\phi\le C\delta_0 |t|^\alpha$ for some $\alpha>0$. We
obtain $|\phi_t|\le C\delta_0|t|^{\alpha-1}$ and so $|y_1|\le \sup
|\phi|\le C\delta_0$.

Next we need a refinement of (3.1).

\proclaim{Lemma 3.3}
 Let $u$ be an entire convex solution of (1.1) with $\inf u=0$. Then
in an appropriate coordinate system, we have (3.1) with
$$|\phi(x)|=O(|x|^{2/3})\ \ \text{as}\ \ |x|\to\infty.\tag 3.13$$
\endproclaim

\noo{\it Proof}.  Let $u_h(x)=h^{-1} u(h^{1/2}x)$. Then $u_h$
satisfies the equation $\Cal L_\sigma [u]=1$ in $\R^2$ with
$\sigma=h^{-1}$. By Theorem 1.3, $u_h$ converges to the function
$u^*= \frac 12|x|^2$, and the level set $\Ga_{1/2, u_h}$ converges
to the unit circle $S^1$ as $h\to\infty$.

For any given sufficiently small constant $\delta_0>0$, let $h>0$
sufficiently large such that
$$ \Gamma_{1/2, u_h} \subset N_{\delta_0}(S^1) \tag 3.14$$
for some unit circle $S^1$. We claim that for any $\tau>\tau_0$,
where $\tau_0>3\max(\delta_0, \sigma)$,
$$\Gamma_{\tau, u_h}\subset
           \sqrt{2\tau}(N_{\delta_\tau}(S^1)) \tag 3.15$$
with
$$\delta_\tau
      \le C_1(\tau) \sigma^{2/3}+C_2\delta_0 \tau^\alpha,$$
where $\alpha (N_\delta(S^1)) = N_{\alpha\delta}(\alpha S^1)$, and
$\alpha S^1$ is the $\alpha$-dilation of $S^1$ with the same center,
the constants $C_1$ and $C_2$ are independent of $\delta_0$ and $h$,
and $C_2$ is also independent of $\tau$. The center of the $S^1$ in
(3.15) is the minimum point of $u_0$, the solution of $\Cal
L_0[u]=1$ in $\Om_{\frac 12, u_h}$ satisfying $u_0=u_h=\frac 12$ on
$\pom_{\frac 12, u_h}$.

To prove (3.15), by Lemma 3.2 we have, for any $\tau>0$,
$$\Gamma_{\tau, u_0}\subset
              \sqrt{2(\tau+a_0)} (N_{\delta_1}(S^1))\tag 3.16$$
with $\delta_1\le C \delta_0 (\tau+a_0)^\alpha$, where $a_0=-\inf
u_0\ge 0$. By the comparison principle we have $u_0\le u_h$ in
$\Om_{\frac 12, u_h}$. By (3.14) we also have
$$u_0\ge \frac 12 (|x|^2-(1+\delta_0)^2)+\frac 12
       \ \ \ \text{in}\ \ \Om_{1/2, u_h}.$$
Hence $a_0\le 3\delta_0$.

We will use the following simple result: Let $\Om$ be a convex
domain contained in $B_R$. If the area $|B_R-\Om|\le \eps$, then
$$\dist (\p B_R, \pom)\le C\eps^{2/3}R^{-1/3},\tag 3.17$$
where $\dist(A, B)$ denotes the least constant $\delta>0$ such
that $A\subset N_\delta (B)$ and $B\subset N_\delta(A)$.

We use (3.17) to prove (3.15). Let $\ell$ be the largest circle,
with center at the minimum point of $u_0$, contained in
$\Om_{\tau, u_0}$. Let $\wtt\Om_{\tau, u}$ be the common area
enclosed by $\Om_{\tau, u}$ and $\ell$, and denote $\wtt\Ga_{\tau,
u}=\p\wtt\Om_{\tau, u}$. Since $\Om_{\tau, u}\subset\Om_{\tau,
u_0}$, we have
$$\dist(\Ga_{\tau, u}, \Ga_{\tau, u_0})\le
   \dist(\wtt\Ga_{\tau, u}, \ell)
    +\dist(\ell, \Ga_{\tau, u_0}).\tag 3.18 $$
Since $\Cal L_0[u_0]=1$, we have $\frac {d}{dt}|\Om_{t,
u_0}|=-2\pi$. Hence $|\Om_{\tau, u_0}|=2\pi (\tau+a_0)$. By
(3.16), $\Ga_{\tau, u_0}$ is in the
$\sqrt{2(\tau+a_0)}\delta_1=C\delta_0 (\tau+a_0)^{1/2+\alpha}$
neighborhood of $\sqrt{2(\tau+a_0)}S^1$. Hence
$$\dist(\ell, \Ga_{\tau, u_0})
        \le C \delta_0 (\tau+a_0)^{1/2+\alpha},$$
where $C>0$ is independent of  $\delta_0$, $h$, and $\tau$. Recall
that  $\Om_{\tau, u}\subset\Om_{\tau, u_0}$ and $|\Om_{\tau,
u_0}-\Om_{\tau, u}|\le C\sigma$ by (3.5). Hence by (3.17) we have
$$\dist(\wtt\Ga_{\tau, u}, \ell)
  \le C \sigma^{2/3}(\tau+a_0)^{-1/6} .$$
Combining (3.16) and (3.18), and noting that $a_0\le
3\delta_0<\tau$, we obtain (3.15).

Now we fix a $\tau_0>0$ small such that $C_2\tau_0^\alpha<1/4$. From
(3.15) we obtain
$$\Gamma_{\tau_0, u_h}
         \subset \sqrt{2\tau_0}(N_{\delta}(S^1))\tag 3.19$$
with $\delta \le C\sigma^{2/3}+\delta_0/4$, where $C$ is
independent of $\delta_0$ and $h$.

Now Lemma 3.3 follows from (3.19) by iteration. We start at the
level $\tau_0^{-k}$ for some sufficiently large $k$. Denote
$\Om_k=\sqrt{2\tau_0^k}\Om_{\tau_0^{-k}, u}$ and $\Ga_k=\p\Om_k$. By
(3.1), $\Ga_k$ converges to the unit circle as $k\to\infty$. Suppose
$\Ga_k$ is in the $\delta_k$-neighborhood of $S^1$, where
$\delta_k\to 0$ as $k\to\infty$. Let $y_k$ denote the minimum point
of the solution of $\Cal L_0[u]=1$ in $\Om_{k+1}$ and $u=\frac 12$
on $\Ga_{k+1}$. By (3.19), $\Ga_{k-1}$ is in the
$\delta_{k-1}$-neighborhood of a unit circle $S^1$ centered at
$y_{k-1}$ with
$$\delta_{k-1}\le C\tau_0^{2(k-1)/3}+\delta_k/4 . $$
By induction, we obtain
$$\delta_{k-2}\le C\tau_0^{2(k-2)/3}+\delta_{k-1}/4.$$
Hence we have
$$\delta_j \le 2C\tau_0^{2j/3}+\delta_k\ \ \ \forall\ j<k.$$
Let $k\to\infty$ we obtain
$$\Ga_j\subset N_{\delta_j}(S^1)\tag 3.20$$
with $\delta_j \le 2C\tau_0^{2j/3}$, where $S^1$ is centered at
$y_j$. It follows that for $h=\tau_0^{-j}$ sufficiently large,
$$\Ga_{h,u}\subset N_{\delta}(\sqrt{2h}S^1)\tag 3.20'$$
with $\delta\le 2Ch^{-1/6}$, where $S^1$ is centered at
$z_j=h^{1/2}y_j$.

Next we estimate $|z_j-z_{j-1}|$. Let
$u_j(x)=\tau_0^ju(\tau_0^{-j/2}x)$. Let $v_1$ and $v_0$ be the
solutions of $\L_0[v]=1$ which satisfy respectively $v_1=1$ on
$\{u_j=1\}$ and $v_0=\tau_0$ on $\{u_j=\tau_0\}$. By (3.20) we have
$\{v_0=\tau_0\}\subset N_\delta(\sqrt{2\tau_0}S^1_{p_0})$ and
$\{v_1=\tau_0\}\subset N_\delta(\sqrt{2\tau_0}S^1_{p_1})$ with
$\delta<C\tau_0^{2j/3}$ for some points $p_0$ and $p_1$. As remarked
before Lemma 3.3, we may assume $p_0$ and $p_1$ are the minimum
points of $v_0$ and $v_1$, so that $z_j=\tau_0^{-j/2}p_1$ and
$z_{j-1}=\tau_0^{-j/2}p_0$. By Lemma 3.1 and (3.17) we also have
$\{v_0=\tau_0\}\subset N_\delta (\{v_1=\tau_0\})$. Hence
$|p_0-p_1|\le C\delta$. We obtain $|z_j-z_{j-1}|\le C\tau_0^{j/6}$.

From the above estimate, the sequence $\{z_j\}$ is convergent.
Assume that $z_j\to 0$. Then the above estimate implies that
$|z_j|\le C\tau_0^{j/6}$ for any large $j$. Hence for
$h=\tau_0^{-j}$,
$$\Ga_{h,u}\subset N_{\delta}(\sqrt{2h}S^1),$$
where $\delta\le Ch^{-1/6}$ and $S^1$ is centered at the origin. It
is easy to see the estimate also holds for all $h>1$. Hence Lemma
3.3 is proved. \hfill$\square$

To finish the proof we need the following fundamental Liouville
Theorem by Bernstein [2], see also [24] (p.245).

\proclaim{Proposition 3.1} Let $u$ be an entire solution to the
elliptic equation
$$\sum _{i, j=1}^2 a_{ij}(x) u_{ij}=0
            \ \ \ \text{in}\ \ \R^2.\tag 3.21$$
If $u$ satisfies the asymptotic estimate
$$|u(x)|=o(|x|)\ \ \ \text{as}\ \ |x|\to\infty, \tag 3.22$$
then $u$ is a constant.
\endproclaim

We remark that the operator in the above proposition need not to
be uniformly elliptic. Condition (3.21) can be replaced by a
weaker condition that $u_{11}u_{22}-u_{12}^2\le 0$ and $\not\equiv
0$.

\noo{\it Proof of Theorem 3.2}. Assume $u$ is locally uniformly
convex, namely the Hessian matrix $(D^2 u)>0$ pointwise, which
will be proved below. Let $u^*$ be the Legendre transform of $u$.
Then $u^*$ satisfies equation (1.9). First we have
$$u^*(x) = \frac 12|x|^2 + O(|x|^{2/3}) . \tag 3.23$$
Indeed, for any $h>1$, let $u_h(x) = h^{-1}u(h^{1/2}x)$. Then by
Lemma 3.3,
$$u_h(x) = \frac 12|x|^2 + O(h^{-2/3})$$
in $B_1(0)$. Denote $u_h^*$ the Legendre transforms of $u_h$. Then
$$u^*_h(x)= \frac 12 |x|^2  +O(h^{-2/3})$$
in $B_1(0)$. Observing that $u^*_h(x)=h^{-1}u^*(h^{1/2}x)$, we
obtain (3.23).

Let $u_0$ be the unique radial solution of (1.1) satisfying
$u(0)=0$, and let $u^*_0$ denote the Legendre transform of $u_0$.
Similar to (3.23) we have
$$u^*_0(x) = \frac 12|x|^2 + O(|x|^{2/3}) . \tag 3.24$$
Write equation (1.9) in the form
$$G[x, D^2 u^*]=:\frac{\det D^2 u^*}
 {\sum(\delta_{ij}-\frac{x_ix_j}{1+|x|^2})F^{ij}[u^*]}=1.\tag 3.25$$
Since both $u^*$ and $u_0^*$ satisfy equation (1.9), $v=u^*-u_0^*$
satisfies equation (3.21) in the entire $\R^2$ with coefficients
$$a_{ij}=\int_0^1 G^{ij}[x, D^2u^*_0+t(D^2u^*-D^2u^*_0)] dt, $$
where $G^{ij}[x, r]=\frac{\p}{\p r_{ij}}G[x, r]$ for any symmetric
matrix $r$. By (3.23) and (3.24), $|v(x)|=O(|x|^{2/3})$ as
$|x|\to\infty$. By the above proposition we conclude that $v$ is a
constant. \hfill$\square$

\vskip5pt

\noo{\bf Remark 3.3}. When using the Legendre transform we have
implicitly used the local uniform convexity of $u$, namely the
Hessian matrix $\{D^2 u\}>0$. In dimension 2, this was proved in
[14] by Hamilton's maximum principle, for high dimensions see [15].
We also note that the reason for using the Legendre transform in the
above proof is that equation (3.21) does not involve the first order
derivatives.

\vskip30pt

%\newpage

\centerline {\bf 4. Translating solutions to the level set flow}

\vskip10pt

In this section we prove the case $\sigma=0$ of Theorem 1.2 and that
an ancient convex (in space) solution to the mean curvature flow is
convex in space-time. We point out that when $n\ge 4$, the proof is
simpler, see Remark 4.1.

\proclaim{Theorem 4.1} For any $n\ge 2$ and $1\le k\le n$, there
exist complete convex solutions, defined in strip regions, to the
equation
$$\sum_{i,j=1}^n (\delta_{ij}-\frac {u_iu_j}{|Du|^2})u_{ij} =1
         \tag 4.1$$
which are not $k$-rotationally symmetric. If $n\ge 3$, there exist
entire convex solutions to (4.1) which are not $k$-rotationally
symmetric.
\endproclaim

By our definition, a function $u$ is $k$-rotationally symmetric if
$u(x)=\phi(|\hat x|)$ in an appropriate coordinate system, where
$\hat x=(x_1, \cdots, x_k)$. To prove Theorem 4.1 we will need the
following logarithm concavity of solutions to (4.1).

\proclaim{Lemma 4.1} Let $\Om$ be a smooth, bounded, convex domain
in $\R^n$. Let $u$ be the solution of (4.1) in $\Om$, vanishing on
$\pom$. Then for any constant $h$ satisfying $\inf_\Om u < h < 0$,
the level set $\Gamma_{h, u}=\{u=h\}$ is convex. Moreover,
$\log(-u)$ is a concave function.
\endproclaim

\noo {\it Proof.}\ \ Since $u$ is a solution of (4.1),
$\psi=-\log(-u)$ satisfies
$$(\delta_{ij}-\frac {\psi_i\psi_j}{|D\psi|^2})\psi_{ij}
               =e^\psi. \tag 4.2$$
Since $\psi(x)\to+\infty$ as $x\to\pom$, the results in [17] (see \S
I\!I\!I.12) implies $\psi$ is convex. \hfill$\square$

Denote
$$\Om_{r, t}=\{x\in\R^n:\ \ \frac{|x'|^2}{r^2}+\frac{x_n^2}{t^2}<1\},
   \tag 4.3$$
where $r, t$ are positive constants, $x'=(x_1, \cdots, x_{n-1})$.
Let $u_{r, t}$ denote the solution of (4.1) in $\Om_{r, t}$,
vanishing on $\pom_{r,t}$. Denote $M_{r, t}=-\inf u_{r, t}$ and
$\Ga_{r, t}=\{u_{r, t}=-M_{r, t}+1\}$. Obviously $M_{r,
t}\to\infty$ as $r, t\to\infty$.

The following lemma plays a key role for our construction of
non-radial convex solutions. A similar idea was used in [4], where
we proved that for any ellipsoid $E$, there exists an entire convex
solution $u$ to the Monge-Amp\`ere equation $\det D^2 u=f$ such that
$u(0)=0$, $u\ge 0$, and the minimum ellipsoid of the sub-level set
$\{u<1\}$ is similar to $E$.

\proclaim{Lemma 4.2} For any $\th>0$ and $K>1$, there exist
$r=r(\th, K)$ and $t=t(\th, K)$ such that $M_{r, t}=K$ and
$$\sup\{|x'|:\ \ x\in \Ga_{r, t}\}
          =\th \sup\{x_n:\ \ x\in \Ga_{r, t}\}. \tag 4.4$$
\endproclaim

\noo{\it Proof}. \ \ The solution $u_{r, t}$ depends continuously
on $r$ and $t$, and $M_{r, t}$ is monotone increasing in $r$ and
$t$. For any $K>1$, we have $M_{r, t}=K$ when
$r=t=\sqrt{2(n-1)K}$.

It is easy to see that for any fixed $r>0$, $M_{r, t}\to 0$ as
$t\to 0$. Hence for any given $r>\sqrt{2(n-1)K}$, there exists a
unique $t=t_r<\sqrt{2(n-1)K}$ such that $M_{r, t}=K$. Moreover we
have $t_r\to 0$ as $r\to\infty$. Similarly for any fixed $t>0$, we
have $M_{r, t}\to 0$ as $r\to 0$. Hence for any given
$t>\sqrt{2(n-1)K}$, there exists a unique $r=r_t<\sqrt{2(n-1)K}$
such that $M_{r, t}=K$.

Observe that for any fixed $K$, $\sup\{x_n:\ \ x\in \Ga_{r, t}\}\to
0$ as $t\to 0$ and by convexity $\sup\{x_n:\ \ x\in \Ga_{r, t}\}\to
\infty$ as $t\to \infty$, where $r=r_t$ is such that $M_{r, t}=K$.
By the continuity of $u_{r, t}$ in $r$ and $t$, there exist $r>0$
and $t=t_r$ such that (4.4) holds. \hfill $\square$

\vskip10pt

For any fixed $\th\ne 1$, by Lemma 4.2 there exist $r=r_k$ and
$t=t_k$ such that $M_{r_k, t_k}=k$ and (4.4) holds. From the proof
of Lemma 4.2 we also have
$$t_k>r_k\ \ \ \ \text{if}\ \ \th<1.\tag 4.5$$
Hence as $k\to\infty$ we have $t_k\to\infty$ if $\th<1$ or
$r_k\to\infty$ if $\th>1$. Denote $w_k=u_{r_k,t_k}+k$. Then
$w_k\ge w_k(0)=0$. We want to prove $w_k$ converges to a complete
solution of (4.1).

We say that a sequence of (embedded) convex hypersurfaces
$\{\M_k\}$ locally converges to $\M$ if for any $R>1$ and
$\delta>0$, there exists $k_0>1$ such that when $k\ge k_0$,
$\M_k\cap B_R(0) \subset N_\delta(\M\cap B_R(0))$ and $\M\cap
B_R(0) \subset N_\delta(\M_k\cap B_R(0))$, where $B_R$ denotes the
ball of radius $R$ and $N_\delta$ denotes the
$\delta$-neighborhood.

For any fixed integer $j$, by Lemma 4.1,  $\phi_{j, k}=\log
j-\log(j-w_k)$ ($k\ge j$) is an even, convex function. Let $\M_{j,
k}$ denotes the graph of $\phi_{j, k}$. Observe that for any fixed
$j$ and $h$, by (4.4) the sets $\M_{j, k}\cap\{x_{n+1}<h\}$ are
uniformly bounded in $k$, and can be represented as radial graphs
with center at the point $(0, \cdots, 0, \frac 12)$. Hence we may
suppose by choosing subsequences that $\M_{j, k}$ converges
locally to a complete, convex hypersurface $\M_j$. Let $\ol D_j$
denote the projection of $\M_j$ on $\{x_{n+1}=0\}$ and $D_j$
denote the interior of $\ol D_j$. Then $D_j$ is a convex domain
and as $k\to\infty$, $\phi_{j, k}$ converges locally in $D_j$ to a
function $\phi_j$. It follows that $w_k$ converges locally in
$D_j$ to a function $w$. Obviously $w$ is a viscosity solution of
(4.1) in $D_j$. Repeating the procedure for $j=1, 2, \cdots$, by
the Arzela-Ascoli lemma we obtain a sequence of domains
$D_1\subset D_2\subset\cdots $ such that $w_k$ sub-converges
locally to $w$ in all $D_j$, $j=1, 2, \cdots$. Let $D=\cup D_j$.
Then $D$ is a convex domain and $w_k$ converges locally to $w$ in
$D$.

By (4.4), $w$ is not rotationally symmetric. To prove Theorem 4.1,
we will prove $w$ is convex and $w(x)\to\infty$ as $x\to\p D$.

\proclaim{Lemma 4.3}\ \ For any $x_0\in\p D$, we have
$$\lim_{x\to x_0} w(x)=\infty. \tag 4.6$$
\endproclaim

\noo {\it Proof}.\ \  For any fixed $k$, the level set $\{w_k=-t\}$
is a convex solution to the mean curvature flow (with time $t\in
(-k, 0)$). For any fixed $t$, by the discussion above we see that
$\{w_k=-t\}$ converges to the level set $\{w=-t\}$ as $k\to\infty$.
Hence $\{w=-t\}$, where $-\infty<t<0$, is also a convex solution to
the mean curvature flow. It follows that for any $t\in(-\infty, 0)$,
$\{w=-t\}$ is smooth and locally uniformly convex. Hence at any time
$t\in (-\infty, 0)$, the hypersurface $\{w=-t\}$ is moving at
positive velocity. Hence $w(x)\to\infty$ as $x\to\p D$.
\hfill$\square$

We can also write the graph of $w_k$ locally in the form
$x_n=v_k(x', t)$, where $t=-x_{n+1}$. Then $v_k$ satisfies the
non-parametric mean curvature flow equation
$$v_t=\sqrt{1+|Dv|^2}\,\div \frac {Dv}{\sqrt{1+|Dv|^2}},\tag 4.7$$
where $Dv=(v_{x_1}, \cdots, v_{x_{n-1}})$. Hence if $v$ is convex
in $x'$ and if $v_t>0$ at some point, then $v_t>0$ everywhere by
the Harnack inequality. Using this property one also easily
conclude (4.6). We remark that Lemmas 4.2 and 4.3 were also
observed by White, see [28].

\proclaim{Lemma 4.4}\ \ The solution $w$ is convex.
\endproclaim

\noo{\it Proof.}\ \ Since the level set of $w_k$ is convex, so is
the level set of $w$. For any point $y\in D$ and any positive
constant $\delta<\min (1, \frac 12 d_y)$, where $d_y=\dist (y, \p
D)$, there exists a constant $M_0>0$ depending on $\delta$ such
that $\sup\{w(x): x\in B_\delta(y)\}\le M_0$ and $\sup\{|D w
(x)|:\ \ x\in B_\delta(y)\}\le M_0$. Denote $v_k=\log(k-w)$.  By
the concavity of $v_k$ we have
$$\align
|Dv_k(y)|
  & \le \sup_{x\in \p B_\delta(y)} \frac 1\delta (|v_k(x)-v_k(y)|)\\
          & \le C(\log (k+M_0) -\log (k-M_0)) \le  C/k , \\
\endalign $$
where $C>0$ depends on $M_0$ and $\delta$, but is independent of
$k$.

By the concavity of $v_k$, we have furthermore
$$\align
\{\p_i\p_j w_k(y)\}
& = -e^{v_k}\{\p_i\p_j v_k +\p_i v_k\p_j v_k\} \\
&\ge -e^{v_k}\{\p_iv_k\p_jv_k\}  \ge - \frac C k I ,\\
\endalign $$
where $I$ is the unit matrix. Sending $k\to\infty$ we obtain
$\{\p_i\p_j w(y)\}\ge 0$. Hence $w$ is convex. \hfill$\square$

From [12] we know that $\{w=-t\}$ shrinks to a round point as $t\to
0$. Hence $w>0$ for any $x\ne 0$ and so $w$ is not $k$-rotationally
symmetric for any $1\le k\le n$. If we choose $\th>1$ sufficiently
large, then $w$ must be defined in a strip region by Lemma 2.7. We
have thus proved the first part of Theorem 4.1.

We would like to point out that, from the proof of Theorems 1.1
and 1.3, the function $w$ is defined in a strip region for any
$\th>1$. If $n=2$, then by Theorem 1.1,  $w$ is defined in a strip
for any $\th\ne 1$.

\vskip10pt

Next we prove the second part of Theorem 4.1.  We will prove the
solution $w$ obtained above is an entire solution if $n\ge 3$ and
$\th<1$. Denote
$$\align
r_h=r_{h,w} & =\sup\{|x'|:\ \ (x', x_n)\in \Om_{h, w}\},\\
t_h=t_{h,w} & =\sup\{x_n:\ \ (x', x_n)\in \Om_{h, w}\}.\\
\endalign $$

\proclaim {Lemma 4.5} Suppose $t_h\ge \delta r_h$ for some
positive constant $\delta>0$. Then
$$\frac {(\delta r_h)^2}{4(n-1)}
   \le h \le \frac {r_h^2}{2(n-2)}.\tag 4.8$$
\endproclaim

\noo{\it Proof}.\ \ Let
$$\phi=\frac {1}{2(n-1)}(|x|^2- \frac 12(\delta r_h)^2).$$
Then $\Cal L_0[\phi]=1$ in $\Om_{h, w}=\{w<h\}$ and $\phi\ge 0$ on
$\p\Om_{h, w}$. By the comparison principle it follows that $w-h
\le \phi$ in $\Om_{h, w}$. Hence
$$h \ge -\inf \phi
        = \frac {(\delta r_h)^2}{4(n-1)}.$$

To prove the second inequality of (4.8), let
$$\phi=\frac {1}{2(n-2)} (|x'|^2-r^2_h).$$
Then $\Cal L_0[\phi]=1$ in $\Om_{h, w}$ and $\phi\le 0$ on
$\pom_{h, w}$. It follows $w-h\ge \phi$ in $\Om_{h, w}$. Hence
$$h\le -\phi(0)= \frac {r_h^2}{2(n-2)}. $$
This completes the proof. \hfill$\square$

Therefore to prove the second part of Theorem 4.1 it suffices to
prove that there exists $\delta>0$ such that
$$\sup\{x_n:\ \ x\in \Ga_{h, w}\}\ge \delta
   \sup\{|x'|:\ \ x\in \Ga_{h, w}\}\tag 4.9$$
for any $h>0$, where $\Ga_{h, w}=\{w=h\}$. Denote
$$\align
r_{h,k} & =\sup\{|x'|:\ \ x\in \Ga_{h, w_k}\},\\
t_{h,k} & =\sup\{x_n:\ \ x\in \Ga_{h, w_k}\}.\\
\endalign $$
Then ${r_{h,k}}_{|\,h=k}=r_k$ and ${t_{h,k}}_{|\,h=k}=t_k$, where
$t_k$ and $r_k$ satisfy (4.5). If there is a subsequence of
$\{k\}$ such that
$$t_{h, k}\ge r_{h, k}
 \ \ \ \forall\ h\in (0, k), \tag 4.10$$
then (4.9) holds with $\delta=1$ for all $h>0$. Hence $w$ is
defined in the entire space $\R^n$.

If (4.10) is not true,  let
$$h_k=\sup \{ h>0:\ \ t_{h, k}<r_{h, k}\}.$$
By (4.5) we have $h_k<k$. If the sequence $\{h_k\}$ is uniformly
bounded, $w$ is defined in the entire space $\R^n$.

If $h_k\to\infty$, we denote $\wtt w_k(x) = h_k^{-1}
w_k(h_k^{1/2}x)$. Then $\Cal L_0 [\wtt w_k]=1$ in $\{\wtt w_k<1\}$
and
$$\sup\{|x'|:\ \ x\in \{\wtt w_k<1\}\}=
   \sup\{x_n:\ \ x\in \{\wtt w_k<1\}\}.\tag 4.11$$
Observe that the level set $\{\wtt w_k=-t\}$ is a convex solution
to the mean curvature flow (with time $t\in (-1, 0)$). From [12]
we see that $\{\wtt w_k=h\}$ shrinks to a round point at $h\to 0$.
Hence by (4.11) we have, for any $h\in (0, 1)$,
$$\sup\{x_n:\ \ x\in \{\wtt w_k=h\}\}\ge \delta
   \sup\{|x'|:\ \ x\in   \{\wtt w_k=h\}\} \tag 4.12$$
for some $\delta>0$ independent of $h$ and $k$. Rescaling back we
obtain (4.9). This completes the proof of Theorem 4.1.

\vskip10pt

\noo{\bf Remark 4.1.} \ When $n\ge 4$, we can also construct entire
convex solutions of (4.1) as follows.  Denote $\hx=(x_1, \cdots,
x_{n-2})$, $\wx=(x_{n-1}, x_n)$. Let $\Om_{r, t}=\{x\in\R^n:\ \
\frac {|\hx|^2}{r^2}+\frac {|\wx|^2}{t^2}=1\}$. Let $u_{r, t}$ be
the solution of (4.1) with $\Om=\Om_{r, t}$, which vanishes on
$\pom$. As before we choose $r_k$ and $t_k$ such that $\inf u_{r_k,
t_k}=-k$ and
$$\sup\{|\hx|:\ \ x\in \Ga_{r, t}\}
          =\th \sup\{|\wx|:\ \ x\in \Ga_{r, t}\}, \tag 4.13$$
where $\Ga_{r, t}=\{u_{r, t}=-k+1\}$, and $\th\ne 1$ is any given
positive constant. Denote $w_k=u_{r_k, t_k}+k$. Then $w_k$ is
nonnegative and $w_k(0)=0$. Observe that the function $\phi=\frac
{1}{2}|\wx|^2$ satisfies the equation $\Cal L_0[\phi]=1$. By the
comparison principle we have
$$w_k(x)\le \frac 12 |x|^2.\tag 4.14$$
Let $w=\lim_{k\to\infty} w_k$. Then $w$ satisfies (4.14). Hence by
Lemma 4.4, $w$ is an entire convex solution of (4.1). Obviously
$w$ is not $k$-rotationally symmetric for any $1\le k\le n$.

\vskip10pt

\noo{\bf Remark 4.2.} \ In Lemma 4.4 we showed that a solution to
(4.1) is convex if the level set $\Ga_h=\{u=h\}$ is bounded and
convex for all large $h>0$. This assertion is true even if $\Ga_h$
is unbounded.

\proclaim{Proposition 4.1} Let $u$ be a solution to (4.1) whose
graph is a complete hypersurface. Suppose that the level set
$\Ga_{h, u}$ is convex for any $h$. Then $u$ is convex.
\endproclaim

\noo{\it Proof}. If the sub-level sets $\Om_{h, u}$ are bounded,
Proposition 4.1 is proved in Lemma 4.4. If $\Om_{h,u}$ is unbounded,
from the proof of Lemma 4.4 it suffices to show that $\log (h-u)$ is
concave for any large constant $h$. By Lemma 4.1, it suffices to
show that for any given $h>0$, $u$ can be approximated locally in
$\Om_{h,u}$ by a sequence of solutions to (4.1) whose level sets are
bounded and convex.

Let $\{D_k\}$ be a sequence of convex domains in $\R^n$ satisfying
$D_1\subset D_2\subset\cdots$, such that $\cup D_k=\Om_{h, u}$ for
some fixed $h$. Let $w_k$ be the solution of (4.1) in $D_k$
satisfying $w_k=h$ on $\p D_k$. Then $w_k\ge u$ in $D_k$ and $w_k$
is decreasing in $k$. Hence $w_k$ converges as $k\to\infty$ to a
solution $w$ of (4.1) in $\Om_h$, satisfying $w=h$ on $\pom_{h,u}$
and $w\ge u$ in $\Om_{h, u}$.

To prove $w=u$, we need to prove that for any $h'<h$ ($h$ fixed),
the level sets $\Ga_{h', w}$ is sufficiently close to $\Ga_{h',
u}$. Note that both $\Ga_{h', w}$ and $\Ga_{h', u}$ evolves by
mean curvature (with time $t=-h'$), so it suffices to prove
$\Ga_{h', w}$ is sufficiently close to $\Ga_{h', u}$ at infinity.

Let $\{x_k\}$ be a sequence on $\Ga_{h, u}$ with $|x_k|\to\infty$.
By the convexity of $\Ga_{h, u}$, the normal of $\Ga_{h, u}$
(regarded as a hypersurface in $\R^n=\{x_{n+1}=h\}$) at $x_k$
converges along a subsequence to a boundary point of the Gauss
mapping image of $\Ga_{h, u}$. Hence after translation, the convex
hypersurface $\Ga^k_{h, u}=\{x-x_k:\ x\in\Ga_{h, u}\}$ converges to
a convex hypersurface which can be split as $\R^1\times \Sigma_{h,
u}$. Similarly for any $h'<h$, $\Ga^k_{h', u}=\{x-x_k:\ x\in\Ga_{h',
u}\}$ converges to a convex hypersurface $\R^1\times \Sigma_{h',
u}$, and $\Ga^k_{h', w}=\{x-x_k:\ x\in\Ga_{h', w}\}$ converges to
$\R^1\times \Sigma_{h', w}$,  and both $\Sigma_{h', u}$ and
$\Sigma_{h', w}$ evolve by mean curvature (with $t=-h'$) with
initial hypersurface $\Sigma_{h, u}$. By an induction argument on
dimension we conclude that $\Sigma_{h', u}=\Sigma_{h', w}$ for any
$h'<h$. Namely $\Ga_{h', w}$ is sufficiently close to $\Ga_{h', u}$
at infinity. \hfill$\square$

\vskip20pt

\centerline {\bf 5. Translating solutions to the mean curvature
flow}

\vskip10pt

In this section we prove the case $\sigma=1$ of Theorem 1.2. That
is

\proclaim{Theorem 5.1} For any dimension $n\ge 2$ and $1\le k\le
n$, there exist complete convex solutions to equation (1.1),
defined in strip regions, which are not $k$-rotationally
symmetric. If $n\ge 3$, there are entire convex solutions to (1.1)
which are not $k$-rotationally symmetric.
\endproclaim

The argument in Section 4 cannot be extended to the mean curvature
equation (1.1), as the logarithm concavity in Lemma 4.1 is still an
open problem for equation (1.1). To prove Theorem 5.1 we will use
the Legendre transform. The purpose to introduce the Legendre
transform is to obtain {\it convex} solutions to the mean curvature
equation (1.1). As remarked at the end of Section 3, we can always
assume that a convex solution is locally uniformly convex.

For clarity we divide this section into three subsections. But
similarly as in \S4, the argument is much more simpler in the case
$n\ge 4$, see discussions at the beginning of \S5.3.

\vskip10pt

\noo{\bf 5.1.} {\it The Legendre transform}

For a smooth, uniformly convex function $u$ defined in a convex
domain $\Om\subset \R^n$. The Legendre transform of $u$, $u^*$, is
a smooth, uniformly convex function defined in the domain
$\Om^*=Du(\Om)$, given by
$$u^*(x)=\sup\{x\cdot y-u(y):\ \ y\in\Om\}. \tag 5.1$$
For example, if $u(x)=a |x|^2$, then $u^*(y)=\frac {1}{4a}|y|^2$;
and if $u(x)=a |x|^{1+\beta}$, then $u^*(y)= c|y|^{1+1/\beta}$ with
$c= {a\beta}/{[a(1+\beta)]^{1+1/\beta}}$. The function $u$ can be
recovered from $u^*$ by the same Legendre transform, namely
$u(y)=\sup\{x\cdot y-u^*(x):\ x\in\Om^*\}$. The supremum is attained
at the unique point $y$ which satisfies
$$x=Du(y) \ \ \ \ \text{and}\ \ \ \ y=Du^*(x).$$
It follows that the Hessian matrix $(D^2u)$ at $y$ is the inverse
of the Hessian matrix $(D^2 u^*)$ at $x$. That is
$$(D^2 u)=(D^2 u^*)^{-1}=(F^{ij}[u^*])/\det D^2 u^*,\tag 5.2$$
where $F^{ij}[u^*]$ is the $(i, j)$-entry of the cofactor matrix
of $(D^2 u^*)$,
$$F^{ij}[u^*]=\frac {\p }{\p r_{ij}}\det\, r
                   \ \ \ \text{at}\ r=D^2 u^*.$$
Hence  if $u$ is a uniformly convex solution of (1.2), $u^*$ is a
solution of $\Cal L^*_\sigma[u^*]=1$, where
$$\Cal L^*_\sigma[u^*]=\det D^2 u^*\big{/}\sum
  (\delta_{ij}-\frac {x_ix_j}{\sigma+|x|^2})F^{ij}[u^*] $$
is a fully nonlinear partial differential equation, which is
elliptic at convex functions.  In particular equation (1.1) is
equivalent to the  equation
$$\Cal L^*_1[u^*]=1.\tag 5.3$$
We have the following classical solvability for the Dirichlet
problem of equation (5.3).

\proclaim{Theorem 5.2} Let $\Om^*$ be a smooth, uniformly convex
domain in $\R^n$ and $\phi$ be a smooth function defined on
$\pom^*$. Then there is a unique, smooth, uniformly convex
solution $u^*\in C^\infty(\bom^*)$ to (5.3) such that $u^*=\phi$
on $\pom^*$.
\endproclaim

For the proof of Theorem 5.2, we observe that the uniqueness of
convex solutions follows from the comparison principle. For the
existence of smooth convex solutions, by the continuity method it
suffices to establish the global regularity estimates. By Evans
and Krylov's elliptic regularity theory, see, e.g., [6,9,18], it
suffices to establish the global second order derivative
estimates.

Different proofs for the global second order derivative estimates
are available [19, 20, 25]. In [19] Krylov provided a
probabilistic proof for (degenerate) Bellman equations. An
analytic proof was later given in [20]. Krylov's estimation covers
equation (5.3) as it is equivalent to a concave equation (5.4)
below and so can be expressed as a Bellman equation. For Hessian
equations (such as (5.4)) the proof in [20] was simplified in
[16].

To apply the a priori estimates in [25] we need to write equation
(5.3) as a Hessian quotient equation, namely
$$\frac  {F_n[w]}{F_{n-1}[w]}=\frac {-1}{p_{n+1}},
              \ \ \ p\in S^*,\tag 5.4$$
where $S^*=\{\frac {(x, -1)}{(1+|x|^2)^{1/2}}\in S^n:\
x\in\Om^*\}$, $S^n$ is the unit sphere,
$$w(p)=(1+|x|^2)^{-1/2}u^*(x)\ \ \ \
   p=\frac {(x, -1)}{(1+|x|^2)^{1/2}}, \ x\in\Om^*.$$
If $u^*$ is the Legendre transform of $u$, the function $w$ is
indeed the support function of $\M_u$ (the graph of $u$), which
can also be defined by
$$w(p)=\sup\{p\cdot X: \ X\in\M_u\}.\tag 5.5$$
Moreover, $\M_u$ can be recovered from $w$ by $\M_u=\p K$ with
$K=\{X\in\R^{n+1}:\ p\cdot X\le w(p)\ \forall\ p\in S^*\}$.

In (5.4) we denote by $F_k[w]$ the $k^{th}$ elementary symmetric
polynomial of the eigenvalues $\lam=(\lam_1, \cdots, \lam_n)$ of
the matrix $\{\D^2w+wI\}(p)$,
$$F_k[w]=\sum_{1\le i_1<\cdots <i_k\le n}
     \lam_{i_1}\cdots\lam_{i_k}, \ \ \ \ 1\le k\le n, \tag 5.6$$
where $\D$ denotes the covariant derivative with respect to a
local orthonormal basis on $S^n$, and $I$ is the unit matrix. As
the supremum in (5.5) is attained at the unique point $X_p\in\M_u$
with normal $p$, the principal radii of $\M_u$ are equal to the
eigenvalues of the Hessian matrix $\{\D^2w+wI\}$ (which also
follows from (5.2)). Hence by equation (1.1), $w$ satisfies
equation (5.4).

The global second order derivative estimates for Hessian quotient
equations in Euclidean domains were established by Trudinger [25].
It is not hard to extend the argument in [25] to equation (5.4)
with domains on the unit sphere.

\vskip10pt

%\newpage

\noo{\bf 5.2.} {\it Complete convex solutions}

With Theorem 5.2 we can now construct a sequence of convex
solutions $(w_k)$ of (1.1), such that $w_k$ converges to a
complete convex solution of (1.1) which is not $k$-rotationally
symmetric for any $1\le k\le n$.

For any positive constants $r, t$, denote
$$\Om^*_{r, t}=\{x\in\R^n:
   \ \ \frac{|x'|^2}{r^2}+\frac{x_n^2}{t^2}<1\} ,$$
where $n\ge 2$, $x'=(x_1, \cdots, x_{n-1})$. By Theorem 5.2, the
Dirichlet problem
$$\cases
\Cal L^*_1[v]=1\ \ &\text{in}\ \ \Om^*_{r, t},\\
 v  =0\ \ &\text{on}\ \ \pom^*_{r,t}\\
 \endcases\tag 5.7 $$
has a unique smooth convex solution $u^*_{r, t}$.  Denote $M^*_{r,
t}=-\inf u^*_{r, t}$ and $\Ga^*_{r, t}=\{x\in\R^n:\ u^*_{r,
t}(x)=-M^*_{r, t}+1\}$. Then $M^*_{r, t}\to\infty$ as $r,
t\to\infty$. Similarly to the proof of Lemma 4.2 we have

\proclaim{Lemma 5.1} For any constants $\th>0$ and $K>1$, there
exist $r=r(\th, K)$ and $t=t(\th, K)$ such that $M^*_{r, t}=K$ and
$$\sup\{|x'|:\ \ x\in \Ga^*_{r, t}\}
       =\th \sup\{x_n:\ \ x\in \Ga^*_{r, t}\}. \tag 5.8$$
\endproclaim

Now we fix a positive constant $\th\ne 1$. By Lemma 5.1 there
exist positive constants $r=r_k$ and $t=t_k$ such that $M^*_{r_k,
t_k}=k$ and (5.8) holds. Similar to (4.5) (after the Legendre
transform the case $\th>1$ here corresponds to the case $\th<1$ in
Section 4) we have
$$r_k>t_k\ \ \ \ \text{if}\ \ \th>1.\tag 5.9$$
Denote $w^*_k=u^*_{r_k,t_k}+k$, $\Om^*_k=\Om^*_{r_k, t_k}$. Then
$w^*_k\ge w^*_k(0)=0$.

Now we use the Legendre transform to change back to equation (1.1).
Let $w_k$ be the Legendre transform of $w^*_k$. Then $w_k$ is a
convex function defined in the domain $\Om_k=:D w^*_k (\Om^*_k)$ and
satisfies the mean curvature equation (1.1) in $\Om_k$. By (5.8),
there exists a constant $R_0>0$, depending only on $n$ and $\th$,
such that $w_k^*\ge 1$ on $\p B_{R_0}(0)$. Hence $|Dw^*_k|\ge 1/R_0$
on $\p B_{R_0}(0)$ and so $B_{R_0^{-1}}(0) \subset \Om_k$ for any
large $k$.

\proclaim{Lemma 5.2} Let $u\in C^2(\Om)$ be a convex solution of
(1.1). Suppose $u(0)=0$, $u\ge 0$, and $u$ is an even function.
Then for any $M>0$, there exists a constant $C>0$ such that for
any $y\in\Om$, if $u(y)<M$, we have
$$|Du(y)|\le C .\tag 5.10$$
\endproclaim

\noo{\it Proof.}\ \
 Write equation (1.1) in the form
$$ \kappa u_\gamma +
 \frac {u_{\gamma\gamma}}{1+u_\gamma^2} =1,\tag 5.11$$
where $\kappa$ is the mean curvature of the level set
$\{u=const\}$ and $\gamma$ is the unit outer normal to the level
set. The normal $\gamma(x)$ is a smooth vector field in
$\Om-\{O\}$. Hence for any point $y\in\Om$, there is a smooth
curve $\ell_y$ connecting the origin $O$ to $y$ such that
$\gamma(x)$ is tangential to the curve at any point $x\in\ell_y$.
Since $u$ is an even function, we may suppose $y$ is in the
positive cone $\{x=(x_1, \cdots, x_n):\ \ x_i\ge 0\}$. It follows
by the convexity of $u$, $\ell_y$ lies in the positive cone and
for any $x\in\ell_y$, $\gamma(x)$ is also a point in the positive
cone. Hence the arc-length $L$ of $\ell_y$ is less than $n|y|$.

Let $\psi$ be the restriction of $u$ on the curve $\ell_y$, and
let $\ell_y$ be parametrized by the arc-length $t$. Then we have
$\psi'=u_\gamma$ and $\psi''=u_{\gamma\gamma}$. Hence
$$\frac {\psi''}{1+{\psi'}^2}=g(t),$$
where $g(t)=1-\kappa \psi'\le 1$. It follows
$$\arctg\, \psi'(t)=G(t)=:\int_0^t g(s)ds, $$
namely $\psi'(t)=\text{tg} G(t)$. Hence $G(L)<\frac \pi 2$ and
$G(L)\ge \frac \pi 4$ if $\psi'(L)\ge 1$. Taking integration we
have
$$\psi(L)=\int_0^L \text{tg} G(t) dt. $$

Denote $L_0=\inf\{t:\ G(t)\ge \frac \pi 4\}$. Then
$$u(y)=\psi(L) \ge \frac {\sqrt 2}{2}
                    \int^L_{L_0} \frac {1}{\cos G(t)}.$$
If $|Du(y)|=\psi'(L)$ is sufficiently large, then $G(L)$ must be
very close to $\frac \pi 2$. It means $\psi(L)$ must be very large
since $G'(t)<g(t)<1$.  Hence Lemma 5.2 is proved. \hfill$\square$

Note that for the relation $\psi''=u_{\gamma\gamma}$ we have used
the fact that $\gamma$ is a normal to the level set $\{u=const\}$.
If $\ell_y$ is replaced by an arbitrary curve $\ell$, then we have
$$\psi''=u_{\gamma\gamma}+\kappa u_\eta,$$
where $\gamma$ is a unit vector tangential to $\ell$, $\kappa$ is
the curvature of $\ell$, and $\eta$ is a unit normal to $\ell$.

Note that by the convexity of $u$, to prove (5.10) it suffices to
consider boundary points. The estimate (5.10) on the boundary of a
convex domain can also be obtained by constructing proper
sub-solutions (more precisely, solution of
$\frac{u_{tt}}{1+u_t^2}=1$ with one variable $t$).

\proclaim{Lemma 5.3} We have
$$m_k=\inf \{w_k(x):\ \ x\in\pom_k\}\to\infty
          \ \ \ \text{as}\ k\to\infty. \tag 5.12$$
\endproclaim

\noo{\it Proof}. First observe that for any constant $h>0$, there
exists $R_h>0$ (independent of $k$) such that for any $x\in\Om_k$
with $|x|>R_h$, we have the estimate
$$w_k(x)\ge h .\tag 5.13$$
In fact this estimate follows from the following geometric
property of the Legendre transform: Let $P_h$ denote the set of
linear functions $g$ such that $g<w_k^*$ and $g(0)=-h$. Let $\ol
g(x)=\sup\{g(x):\ g\in P_h\}$. Then the graph of $\ol g$ is a
convex cone and $D\ol g(\R^n) = \{w_k\le h\}$, where $D\ol
g(\R^n)$ is the image of the sub-gradient mapping,
$$D\ol g(\R^n)=D\ol g(\{0\})
 =\{p\in\R^n:\ \ol g(x)\ge p\cdot x +\ol g(0)
 \ \forall\ x\in\R^n\}. $$
By (5.8) we have $|D\ol g|<C$. Namely for any $h>0$, the set
$\{w_k<h\}$ is uniformly bounded. Hence (5.13) holds.

It follows by convexity that for any boundary point
$x_k\in\pom_k$,  we have $w_k(x_k)\to\infty$ if $|x_k|\to \infty$.
Therefore to prove (5.12) we need only to consider any bounded
sequence $x_k\in\pom_k$.

However, if $|x_k|$ and $w_k(x_k)$ are both uniformly bounded,
then by Lemma 5.2, $y_k=Dw_k(x_k)\in\pom^*_k$ are also uniformly
bounded. Hence by the Legendre transform,
$$k=w^*_k(y_k) = x_k\cdot y_k-w_k(x_k)$$
are also uniformly bounded, a contradiction. \hfill$\square$

\proclaim{Lemma 5.4} There is a subsequence of $\{w_k\}$ which
converges to a complete convex solution $w$ of (1.1).
\endproclaim

\noo{\it Proof}. For any constant $h>0$, by (5.13) the sets
$\M_{w_k}\cap\{x_{n+1}<h\}$ are uniformly bounded in $k$, where
$\M_{w_k}$ denotes the graph of $w_k$. Hence we may suppose by
choosing subsequences that $M_{w_k}$ converges locally to a convex
hypersurface $\M$, which is complete by Lemmas 5.2 and 5.3. Let $\ol
D$ be the projection of $\M$ on the plane $\{x_{n+1}=0\}$ and $D$
the interior of $\ol D$. Then $D$ is a convex domain and $w_k$
converges locally in $D$ to a convex function $w$, and $w$ is a
convex solution of (1.1) in $D$. We claim that $w(x)\to\infty$ as
$x\to\p D$. Indeed, if this is not true, then there is a point $p\in
\M$ at which the tangent plane of $\M$ is perpendicular to the plane
$\{x_{n+1}=0\}$. Hence there is a sequence $p_k=(x_k,
w_k(x_k))\in\M_{w_k}$, $p_k\to p$, such that $w_k(x_k)$ is uniformly
bounded but $|Dw_k(x_k)|\to\infty$. This is in contradiction with
Lemma 5.2.  Consequently $w$ is a complete convex solution of (1.1).
By (5.8), $w$ is not rotationally symmetric. \hfill$\square$

Now we can prove the first part of Theorem 5.1.  Indeed, let $P$
denote the set of linear functions $g$ such that $g<w$ and
$g(0)=-1$. Let $\ol g(x)=\sup\{g(x):\ g\in P\}$. Then the graph of
$\ol g$ is a convex cone and $D\ol g(\R^n) = \{w^*\le 1\}$, where
$w^*=\lim_{k\to\infty} w^*_k$ is the Legendre transform of $w$.
Hence if the constant $\th>0$ in Lemma 5.1 is chosen sufficiently
small, then the level set $\Ga_{\ol g}=\{x\in\R^n:\ \ol g(x)=1\}$
satisfies
$$\sup\{|x'|:\ x\in\Ga_{\ol g}\}
       =\th' \sup\{x_n:\ x\in\Ga_{\ol g}\}$$
for some $\th'>1$ sufficiently large. Hence by Lemma 2.7, $w$ is
defined in a strip region.

\vskip10pt

\noo{\bf 5.3.} {\it Entire convex solutions}

Next we prove the second part of Theorem 5.1. We prove that if
$n\ge 3$ and $\th>1$, the solution $w$ obtained in Lemma 5.4 is
defined in the entire space $\R^n$. The following proof is
necessary only when $n=3$, since if $n\ge 4$, one can construct a
sequence of functions $w^*_k$ as above such that $w_k$, the
Legendre transform of $w_k^*$, takes the form $w_k(x)=w_k(|\hat
x|, |\wtt x|)$ (where $\hat x=(x_1, \cdots, x_{n-2})$, $\wtt
x=(x_{n-1}, x_n)$, see Remark 4.1). Then $w_k(x)\le \frac 12
|x|^2$ and so $\{w_k\}$ sub-converges to an entire convex solution
of (1.1).

For any $h>0$, denote
$$ \align
a_{h,k} & =\sup\{|x'|:\ \ x \in \Ga_{h, k}\} ,\tag 5.14\\
b_{h,k} & =\sup\{x_n:\ \ x\in \Ga_{h, k}\}, \\
\endalign $$
where $\Ga_{h,k}=\{w_k=h\}$.  Let
$$\hat w=\frac {1}{2(n-2)}(|x'|^2-a_{h,k}^2)+h, $$
we have $\Cal L_1[\hat w]\ge 1$. By the comparison principle we
have $w_k\ge \hat w$ and so $\hat w(0)\le 0$. We obtain $a_{h,
k}\ge (2(n-2)h)^{1/2}$. Sending $k\to\infty$ we obtain
$$a_{h}\ge (2(n-2)h)^{1/2},\tag 5.15$$
where we denote
$$ \align
a_h & =\sup\{|x'|:\ \ x \in \Ga_{h}\} ,\\
b_h & =\sup\{x_n:\ \ x\in \Ga_{h}\}, \\
\endalign $$
where $\Ga_h=\{w=h\}$.  To prove that $w$ is an entire solution,
it suffices to prove
$$b_h\to\infty \ \ \text{as}\ \ h\to\infty .\tag 5.16$$

To prove (5.16) we use Lemma 2.7. That is if there exist $h_0>1$ and
$\beta>0$ sufficiently small,  such that $b_{h_0,k}\le \beta
h_0^{1/2}$, then there exists a constant $C>0$ independent of $k$,
such that $b_{h,k}\le C$ for all $h\in (1, m_k)$, namely $w_k$ is
defined in the strip $\{|x_n|<C\}$, where $m_k$ is the constant in
(5.12).

It follows that if $b_{\hat h,k}\ge \beta \hat h^{1/2}$ for some
large $\hat h>1$, then there exists $\delta>0$, independent of
$\hat h$ and $k$, such that $b_{h,k}\ge \delta h^{1/2}$ for all
$1<h<\hat h$. Therefore to prove $w$ is an entire solution, it
suffices to prove that there exists a constant $\beta>0$
independent of $k$, and a sequence $\tau_k$, where $\tau_k\le m_k$
and $\tau_k\to\infty$ as $k\to\infty$, such that
$$b_{\tau_k, k}\ge \beta \tau_k^{1/2}.\tag 5.17$$

Denote
$$ \align
r_{h,k} & =\sup\{|x'|:\ \ x \in \Ga^*_{h, k}\} ,\tag 5.18\\
t_{h,k} & =\sup\{x_n:\ \ x\in \Ga^*_{h, k}\}, \\
\endalign $$
where $\Ga^*_{h,k}=\{w^*_k=h\}$. We have

\proclaim{Lemma 5.5} Suppose $r_{h, k}\ge t_{h,k}$. Then we have
the estimates
$$\sqrt{h/n}\le r_{h, k}\le \sqrt {2h}. \tag 5.19$$
\endproclaim

\noo{\it Proof}.\ The function
$$v=h+\frac n2(|x|^2-2r^2_{h, k})$$
is a sub-solution of (5.3), namely $\Cal L^*_1[v]\ge 1$. Moreover,
since $r_{h, k}\ge t_{h, k}$,  we have $v\le h$ on $\{w^*_k=h\}$.
Hence by the comparison principle, $v\le w^*_k$. In particular we
have $0=w^*_k(0)\ge v(0)= h-nr_{h, k}^2$. The first inequality of
(5.19) is proved.

Next observe that when $n\ge 3$, the function
$$v=h+\frac 12(|x'|^2-r_{h, k}^2)+ K x_n^2$$
is a super-solution of (5.3), namely $\Cal L^*_1[v]\le 1$, for any
$K>1$. For any $\eps>0$ we can choose $K$ sufficiently large such
that $v>w_k^*-\eps$ on $\{w_k^*=h\}$. Hence $v(0)\ge w^*_k(0)=0$ and
we obtain the second inequality. \hfill$\square$

We remark that $\L_1^*[u^*]\ge 1$ ($\le 1$, resp.) is equivalent
to $\L_1[u]\le 1$ ($\ge 1$, resp.), where $u^*$ is the Legendre
transform of $u$.

\vskip10pt

\noo{\it Proof of Theorem 5.1}. The first part of Theorem 5.1 has
been proved in \S 5.2. We need only to prove that when $n\ge 3$
and $\th>1$, $w$ is an entire solution. It suffices to prove
(5.17).

By (5.9) we have $r_{h, k}\ge t_{h, k}$ when $h$ is close to $k$.
If $r_{h, k}\ge t_{h, k}$ for all $h<k$, then by (5.19) we have
$$\frac 12 |x'|^2 \le w^*_k(x', 0)\le n |x'|^2. \tag 5.20$$
By the Legendre transform we have
$$\frac 1{4n} |x'|^2 \le w_k(x', 0)\le \frac 12 |x'|^2. $$
Sending $k\to\infty$ we obtain
$$\frac 1{4n} |x'|^2 \le w(x', 0)\le \frac 12 |x'|^2. \tag 5.21$$
Hence
$$a_h\le (4nh)^{1/2}\ \ \ \forall\ h>1. $$
By Lemma 2.6,  $a_hb_h\ge C h$ ($a_h$ and $b_h$ here correspond to
$\ol a_h$ and $\ol b_h$ in Lemma 2.6). Hence $b_h\to\infty$ as
$h\to\infty$, namely $w$ is defined in the whole space $\R^n$.

If there exists $h>0$ such that  $r_{h, k}< t_{h, k}$, we denote
$$h_k=\inf\{h:\ \ r_{h', k}\ge t_{h', k}
                    \ \ \forall\  h<h'<k\}.\tag 5.22$$
If $\{h_k\}$ is uniformly bounded, or contains a uniformly bounded
subsequence, then (5.20), and so also (5.21), holds for $|x'|$
large. Hence $w$ is also an entire solution.

Finally we consider the case $h_k\to\infty$ as $k\to\infty$.
Denote $G_k^*=\{w^*_k<h_k\}$ and $G_k=Dw^*_k(G^*_k)$. We consider
the solution $w_k$ in the domain $G_k$. Denote
$$\tau_k=:\inf\{w_k(x):\ x\in  \p G_k\}. \tag 5.23$$
Since $h_k\to\infty$, similarly to Lemma 5.3 we have $\tau_k\to
\infty$ as $k\to\infty$. Denote $\hat G^*_k=\frac {1}{\sqrt {h_k}}
G^*_k$, $\hat G_k=\frac {1}{\sqrt {h_k}} G_k$, and
$$\align
\hat w^*_k(x) & =\frac {1}{h_k} w^*_k(\sqrt{h_k} x),\\
 \hat w_k(x) & =\frac {1}{h_k} w_k(\sqrt{h_k} x). \tag 5.24\\
 \endalign $$
Then $\hat w_k$ is the Legendre transform of $\hat w^*_k$, and
$$c_k=:\inf\{\hat w_k(x):\ x\in \p\hat G_k\}=\tau_k/h_k.\tag 5.25$$
We claim
$$c_k\le 4 n \ \ \ \forall\ k.\tag 5.26$$
Indeed, by (5.22), $\hat G^*_k$ has a good shape, namely
$$\sup\{|x'|:\ x\in \hat G^*_k\}
           =\sup\{x_n:\ x\in \hat G^*_k\}.\tag 5.27$$
As the domain $\hat G^*_k$ is rotationally symmetric with respect
to $x'$, by Lemma 5.5 and (5.27) we have $\hat G^*_k\subset
B_2(0)$. Note that $\hat w^*_k$ is the Legendre transform of $\hat
w_k$, we have $D\hat w_k(\hat G_k)=\hat G^*_k$. Hence
$$|D \hat w_k|\le 2\ \ \ \text{in}\ \ \hat G_k.\tag 5.28$$
Since $\hat G^*_k\subset B_2(0)$, there exists $r\le 2$ such that
$\hat G^*_k\subset B_r(0)$, and $\p \hat G^*_k$ and $\p B_r(0)$ have
a common boundary point $x_0$. Observe that $v=1+\frac
n2(|x|^2-r^2)$ is a sub-solution of (5.3) and $v\le \hat w^*_k=1$ on
$\p\hat G^*_k$, we have $|D\hat w^*_k(x_0)|\le |Dv(x_0)|\le 2n$. By
the Legendre transform, $y_0=D\hat w^*_k(x_0)$ is a boundary point
of $\hat G_k$. Since $\hat w^*_k$ is constant on $\p \hat G^*_k$,
one easily verifies that the domain $\hat G_k=D\hat w^*_k(\hat
G^*_k)$ is star-shaped. By (5.28) we obtain
$$\hat w_k(y_0)\le |y_0|\sup |D\hat w_k|\le 4n,$$
namely (5.26) holds.

Denote
$$\align
\eps_k & =\inf\{x_n>0:\ \ \hat w_k(0, \cdots, 0, x_n)\ge c_k\}\\
   & =\sup\{x_n:\ \ \hat w_k(x', x_n)\le c_k\} . \\
   \endalign $$
Then by (5.26) and the rescaling (5.24),  to prove (5.17) it
suffices to prove that $\eps_k$ has a uniform positive lower
bound, namely
$$\eps_k  \ge \eps_0>0\ \ \ \forall\ k \tag 5.29 $$

We prove (5.29) as follows. Observing that $\hat w^*_k(0)=0$,
$\hat w^*_k=1$ on $\p \hat G^*_k$, and $\hat G^*_k\subset B_2(0)$,
we have $|D\hat w^*_k|\ge 1/2$ on $\p\hat G^*_k$. That is $\hat
G_k\supset B_{1/2}(0)$ for all $k$. By (5.28) we may suppose $\hat
w_k\to \hat w$ in $B_{1/2}(0)$. Then $\hat w\ge 0$, $\hat w(0)=0$,
and $|D\hat w|\le 2$ in $B_{1/2}(0)$. Since $\hat w_k$ satisfies
the equation $\Cal L_\sigma[v]=1$ in $\hat G_k$, where
$\sigma=1/h_k$, and $\Cal L_\sigma$ is the operator in (1.2). By
Lemma 2.5, $\hat w$ is a solution of
$$\Cal L_0[v]=1 \ \ \text{in}\ \ B_{1/2}(0).\tag 5.30 $$

If $\hat w(x)>0$ for any $x\ne 0$, then $c=\inf\{\hat w(x):\
|x|=1/2\}>0$. Hence $c_k\ge \frac 12 c$ for sufficiently large
$k$. By (5.28) we then have $\eps_k\ge c_k/2\ge c/4$ for large
$k$. Hence (5.29) holds.

If the convex set $\{\hat w=0\}$ contains a line segment $\ell$,
then $\ell$ is either contained in the $x_n$-axis, or in the plane
$\{x_n=0\}$, as the function $\hat w$ is rotationally symmetric
with respect to $x'$ and symmetric in $x_n$. Furthermore, by the
comparison principe, the set $\{\hat w=0\}$ contains no interior
points, i.e. the Lebesgue measure $|\{\hat w=0\}|=0$. Note that
the origin is the middle point of $\ell$ as $\hat w$ is an even
function.

In the former case, namely if $\ell$ is contained in the
$x_n$-axis, we have $\hat w(x)> 0$ for any point $x\in\p
B_{1/2}(0)$ not lying on the $x_n$-axis, for otherwise $\hat w=0$
on a set with interior points. Hence
$$c_k  \ge \inf_{\p Q}\hat w_k(x)
       = \hat w_k(0, \cdots, 0,\frac d2) , $$
where $Q=\{x\in\R^n:\ |x'|<\frac d2, |x_n|<\frac d2\}$, where $d$
is the arclength of $\ell$. Hence we have $\eps_k\ge \frac d2$ and
(5.29) holds.

We claim the latter case, namely the case when the line segment
$\ell$ is contained in the plane $\{x_n=0\}$, does not occur.
Indeed, if $\ell\subset\{x_n=0\}$, then $\{\hat w=0\}$ is a disc
type set $D=\{x\in B_{1/2}(0):\ |x'|<d, x_n=0\}$ for some $d>0$.
Since the set $\{\hat w=0\}$ contains no interior points, the level
set $\{x\in\R^n:\ \hat w(x)= h\}$ is contained in a strip
$\{|x_n|<\delta\}$ with $\delta\to 0$ as $h\to 0$. Therefore if
$h>0$ is sufficiently small, there exist points on the level set
$\{\hat w= h\}$ at which the mean curvature $\kappa$ of $\{\hat
w=h\}$ is sufficiently small. Write equation (5.30) in the form
$\kappa\, |D\hat w|=1$ (see (2.2) with $\sigma =0$), we find that
$|D\hat w|$ is very large. On the other hand we have $|D\hat w|\le
2$ by (5.28). We reach a contradiction. Hence the latter case does
not occur. This completes the proof. \hfill$\square$

\vskip20pt

\baselineskip=13.6pt
\parskip=2.5pt

\centerline{\bf 6. Applications to the mean curvature flow}

\vskip10pt

First we consider applications of Theorem 1 1. For a mean convex
flow $\M=\cup_{t\in [0, T)}\M_t$ in $\R^{n+1}$, let $\Cal F$ denote
the set of all limit flows (namely blowup solutions) to $\M$ before
first time singularity. A key result in [27] is that a limit flow in
$\Cal F$ cannot be a hyperplane of multiplicity two, from which it
follows that the $``$grim reaper$"$ $x_{n+1}= \log\sec x_1$ can not
be a limit flow in $\Cal F$, see [28, 30].  As indicated in the
introduction, we will consider solutions with positive mean
curvature only.

\proclaim{Corollary 6.1} A flow $\M'\in \Cal F$ (with positive
mean curvature at some point) must sweep the whole space
$\R^{n+1}$.
\endproclaim

\noo{\it Proof}. By Proposition 4.1, $\M'$ is a graph in space-time
(with $x_{n+2}=-t$) of a convex function $u$ on a convex domain in
$\R^{n+1}$, and $u$ is a complete convex solution of (1.2) with
$\sigma=0$. If $u$ is not defined in the entire $\R^{n+1}$, by
Corollaries 2.1 and 2.2, $u$ is defined in a strip region of the
form $\{x\in\R^{n+1}:\ |x_{n+1}|< C\}$ (in appropriate coordinates).
By Lemma 2.6, the projection of $\{u=h\}$ on the $x'=(x_1, \cdots,
x_n)$ plane contains the ball $\{|x'|<Ch\}$ for some $C>0$
independent of $h\ge 1$. It follows that the tangent flow at
infinity to the solution $u$ (more precisely, there is a limit flow
to the original mean convex flow which) is a multiplicity 2 plane,
see Corollary 12.5 in [27] or Corollary 4 in [28]. But a
multiplicity 2 plane does not occur as a limit flow in $\Cal F$.
Hence $u$ must be defined in the whole $\R^{n+1}$. \hfill$\square$

If a convex translating solution $\M'$ is a limit flow to a mean
convex flow in $\R^3$, then $\M'$ is the graph of a convex
function $u$ satisfying (1.1). By Corollary 6.1, $u$ is defined in
the whole $\R^2$. By Theorem 1.1, $u$ is rotationally symmetric.
We obtain Corollary 1.1.

From the case $\sigma=0$ in Theorem 1.1, we also have the
following

\proclaim{Corollary 6.2}  A convex solution to the curve
shortening flow which sweeps the whole space $\R^2$ must be a
shrinking circle.
\endproclaim

From Theorem 1.2, there exist closed convex solutions to the curve
shortening flow which are not the shrinking circle.

\vskip10pt

Next we consider applications of Theorem 1.3. First we have

\proclaim{Corollary 6.3}  Let $\M=\{\M_t\}$ be an ancient convex
solution to the mean curvature flow. Let $\M_t'=\{x\in R^{n+1}:\
(-t)^{1/2}x\in \M_t\}$ be a dilation of $\M_t$. Then $\M'_t$,
after a proper rotation of coordinates, converges as $t\to-\infty$
to one of the following
\newline
(i) an $n$-sphere of radius $\sqrt {2n}$;
\newline
(ii) a cylinder $S^k\times \R^{n-k}$, where $S^k$ is a $k$-sphere
of radius $\sqrt{2k}$;
\newline
(iii) the plane $\R^n$ of multiplicity two.
\endproclaim

\noo{\it Proof}. By Proposition 4.1,  $\M$ is the graph in
space-time of a convex function $u$ in $\R^{n+1}$ satisfying
equation (1.2) with $\sigma=0$. If $u$ is an entire solution, (iii)
does not occur and Corollary 6.3 is equivalent to Theorem 1.3. If
$u$ is not an entire solution, by Corollaries 2.1 and 2.2, $u$ is
defined in a strip region. As in the proof of Corollary 6.1, the
projection of $\M_t=\{u=-t\}$ on the $x'$-plane contains the ball
$\{|x'|<C|t|\}$, by convexity we see that $\M'_t$ converges to a
plane. \hfill$\square$

From Theorem 1.3 (the case $\sigma=0$) we also obtain

\proclaim{Corollary 6.4}  Let $\M=\bigcup_{t\in[0, T)} \M_t$ be a
mean convex flow in $\R^{n+1}$. Suppose $h_i\to\infty$ and
$p_{i}\in M_{t_{i}}$ such that the blow-up sequence
$\M_i=\{(h_i(p-p_{i}), h_i^2(t-t_{i})):\ (p, t)\in \M\}$ converges
to an ancient convex solution. Then there exist $h'_i\to\infty$
such that the corresponding blow-up sequence $\M_i'=\{
(h'_i(p-p_{i}), {h'}_i^2(t-t_{i})):\ (p, t)\in \M\}$ converges
along a subsequence to a shrinking sphere or cylinder.
\endproclaim

At type I\!I singularity, from Theorem 1.3 (the case $\sigma=1$)
we have

\proclaim{Corollary 6.5}  Let $\M_i$ be as in Corollary 6.4. If
$\M_i$ converges to a convex translating solution, then there is a
sequence of positive constants $\lam_i\to 0$ such that, in a
proper coordinate system, the blow up sequence
 $\wtt\M_i=\{ (\lam_i x_1, \cdots, \lam_ix_{n},
    \lam_i^2x_{n+1},\lam_i^2 t): \ (x_1, \cdots, x_{n+1}, t)\in \M_i\}$
converges to the flow $\M'=\{(x, \eta_k(x)+t)\in\R^{n+1}:\ x\in
\R^n, t\in\R\}$, where $\eta_k$ is given in (1.4).
\endproclaim

In Corollaries 6.3-6.5, we need not to restrict ourself to limit
flows in $\F$.

As indicated in the Introduction, our Theorem 1.3 corresponds to
Perelman's classification of ancient $\kappa$-noncollapsing
solutions with nonnegative sectional curvature to the 3-dimensional
Ricci flow [22]. Indeed, Theorems 2.1 and 2.2 implies that the set
of entire, ancient convex solutions in $\R^{n+1}$, which includes
all limit flows in $\F$ by Corollary 6.1, is compact if one
normalizes the solutions such that their mean curvature equals one
at some fixed point, say the origin. To see this, let $\M=\{\M_t\}$
be an ancient convex solution. Then $\M$ is the graph in space-time
of an entire convex function $u$. If the mean curvature of $\M$
equals one at the origin, then $|Du(0)|=1$, and so the compactness
follows from Theorems 2.1 and 2.2.

For ancient convex solutions in $\R^3$, we have the following

\proclaim{Corollary 6.6} Let $\M=\{\M_t\}$ be an ancient convex
solutions which is a limit flow to a mean convex flow in $\R^3$.
Suppose at time $t=0$, the time slice $\M_0=\{u=0\}$ is noncompact.
Then $\forall$  $\eps>0$, there is a compact set $G$ such that any
point in $\M_0\back G$ has a neighborhood which is, after
normalization, $\eps$-close (see (2.56)) to the cylinder $S^1\times
\R^1$.
\endproclaim

\noo{\it Proof}. As indicated above, $\M$ is the graph in space-time
of a convex function $u$ defined in the whole $\R^3$. If the level
sets $\{u=h\}$ is a cylinder, then $u$ is a convex function of two
variables and Corollary 6.6 follows from Theorem 1.1. In this case
the set $G=\emptyset$. Otherwise by convexity we may assume by
choosing an appropriate coordinate system that $u(0)=0$, $u\ge 0$ in
$\{x_1\le 0\}$, and $u\le 0$ on the positive $x_1$-axis. For $a\in
(0, \infty)$, let $x_a$ be the point such that $u(x_a)=\inf \{u(x):\
x_1=a\}$. If $\inf u$ is bounded, let $u_a(x)=u(x+x_a)-u(x_a)$. By
the convexity of $u$, $|Du(x_a)|$ is decreasing as $a\to\infty$, so
it is uniformly bounded. By the above mentioned compactness result
(Theorem 2.1), $u_a$ sub-converges to a convex function $u_0$. By
our choice of coordinates, we have $u_0\le 0$ on the positive
$x_1$-axis. But since $\inf u$ is bounded, we have $u_0= 0$ on the
positive $x_1$-axis. By Lemma 2.9 it follows that $u_0$ is
independent of $x_1$. That is $u_0$ is a function of $x_2$ and
$x_3$. By Theorem 1.1, $u_0$ is rotationally symmetric in $x_2$ and
$x_3$.

If $\inf u=-\infty$, then $u(x_a)\to-\infty$ as $a\to\infty$. By the
compactness Theorem 2.1, the sequence $u_a(x)=\frac 1{h_a}[u(\sqrt
{h_a}x+x_a)-u(x_a)]$, where $h_a=|u(x_a)|$, converges along a
subsequence to a convex solution $u_0$ of (4.1). Since $|Du(x_a)|$
is uniformly bounded, we see that $u_0=0$ on the positive
$x_1$-axis, see (2.49). By Lemma 2.9 it follows that $u_0$ is
independent of $x_1$. By Theorem 1.1, $u_0$ is rotationally
symmetric in $x_2$ and $x_3$. \hfill$\square$

By a compactness argument, one easily sees that the diameter of $G$
is uniformly bounded if the maximum of the mean curvature of $\M_0$
is equal to 1. We remark that by Theorem 1.2, Corollary 6.6 is not
true for convex solutions in $\R^n$ for $n\ge 4$. But in high
dimensions we have accordingly Corollary 6.3, which says that an
ancient convex solution behaves asymptotically as $t\to-\infty$ like
a sphere or cylinder. Note that for any dimension $n\ge 2$, the set
$\F$ contains all blowup solutions before the first time
singularity, and that the blowup sequence converges smoothly on any
compact sets to an ancient convex solution [28]. Therefore similarly
as in [23], one may infer that if $\M=\{\M_t\}$ is a mean convex
flow, then at any point $x_t\in\M_t$ with large mean curvature
before the first time singularity, $\M_t$ satisfies a canonical
neighborhood condition. If $n=2$, the condition is very similar to
that in [23]. In high dimensions, the condition is more complicated.
We will not get into details in this direction. Concerning the
geometry of singularity set it would be more interesting if one can
prove the following

\proclaim{Conjecture} Let $K$ be a smooth and compact domain in
$\R^{n+1}$. Suppose the boundary $\p K$ has positive mean
curvature. Let $\M=\cup_t \M_t$ be the solution to the mean
curvature flow with initial condition $\p K$. Then
\newline
(i) singularity occurs only at finitely many times;
\newline
(ii) at each singular time the singularity set consists of
finitely many connected components;
\newline
(iii) each connected component is contained in a $C^1$ smooth
$(n-1)$-submanifold (with or without boundary). If $n=2$, each
component is either a single point or a $C^1$ curve.
\endproclaim

One may wish to assume $n<7$, but the conjecture is likely to be
true for all $n\ge 2$, as the second order derivative estimates for
(4.1) may hold for any $n\ge 2$.

\vskip10pt

We conclude this paper with some interesting questions related to
our Theorems 1.1-1.3. For Theorem 1.1 a question is whether an
entire solution of (1.1) in $\R^2$ is convex. For Theorem 1.2 a
question is whether a convex solution $u$ of (1.1) is rotationally
symmetric if $|Du(x)|\to\infty$ as $|x|\to\infty$. We expect
affirmative answers to both questions.

For Theorem 1.3, an interesting question is whether a
non-rotationally symmetric ancient convex solution can occur as a
limit flow in $\R^n$. We believe that any limit flow to a mean
convex flow at isolated singularities in space-time is
rotationally symmetric, otherwise non-rotationally symmetric
convex limit flow may occur if the following situations arise: (a)
if there exists a mean convex flow in $\R^n$ ($n\ge 4$) which
develops first time type I\!I singularities simultaneously on a
non-smooth curve (say a polygon) in a 2-plane; (b) if a mean
curvature flow in $\R^n$ ($n\ge 4$) develops first time
singularity on a smooth curve, and the singularity is type I,
except one type I\!I singular point. In case (a) we expect a
non-rotationally symmetric convex translating solution at the
vertices of the polygon. In case (b) a blowup solution near the
type I\!I singular point may not be rotationally symmetric.

 \vskip30pt

\parskip=0.5pt

\Refs\widestnumber\key{123}

\item {[1]}   B. Andrews,
              Contraction of convex hypersurfaces in Euclidean space,
              Calc. Var. PDE, 2 (1994), 151--171.

\item {[2]}   S.N. Bernstein,
              Sur un theoreme de geometrie et ses applications
              aux equations aux derivees partielles du type
              elliptique,
              Comm. de la Soc. Math de Kharkov (2eme ser.),
              15(1915-17), 38-45.
              See also:
              \"Uber ein geometrisches Theorem und seine Anwendung
              auf die partiellen Differential gleichungen vom
              elliptischen Typus,
              Math. Z. 26(1927), 551-558.

\item {[3]}   Y.G. Chen, Y. Giga,  and S. Goto,
              Uniqueness and existence of viscosity solutions
              of generalized mean curvature flow equations,
              J. Diff. Geom., 33(1991), 749-786.

\item {[4]}   K.S. Chou and X.-J. Wang,
              Entire solutions to the Monge-Amp\`ere equations,
              Comm. Pure Appl. Math., 49(1996), 529-539.

\item {[5]}   K. Ecker and G. Huisken,
              Mean curvature evolution of entire graphs,
              Ann. of Math., 130 (1989), 453--471.

\item {[6]}   L. C. Evans,
              Classical solutions of fully nonlinear,
              convex, second-order elliptic equations.
              Comm. Pure Appl. Math. 35 (1982), 333--363.

\item {[7]}   L. C. Evans and J. Spruck,
              Motion of level sets by mean curvatures,
              J. Diff. Geom., 33(1991), 635-681.

\item {[8]}   M. Gage and R.S. Hamilton,
              The heat equation shrinking convex plane curves.
              J. Diff. Geom. 23(1986), 69--96.

\item {[9]}  G. Gilbarg and N.S. Trudinger,
             Elliptic partial differential equations of
             second order,
             Second edition, Springer-Verlag, 1983.

\item {[10]} R.S. Hamilton, Formation of singularities in the
             Ricci flow, Surveys in Diff. Geom., 2(1995), 7-136.

\item {[11]} R.S. Hamilton,
             Harnack estimate for the mean curvature flow.
             J. Diff. Geom., 41(1995), 215--226.

\item {[12]} G. Huisken,
             Flow by mean curvature of convex surfaces into spheres.
             J. Diff. Geom., 20(1984), 237--266.

\item {[13]} G. Huisken,
             Local and global behavior of hypersurfaces moving by
             mean curvature.
             Proc. Symp. Pure Math., 54(1993), 175--191.

\item {[14]} G. Huisken and C. Sinestrari,
             Mean curvature flow singularities for mean convex surfaces,
             Calc. Var. PDE, 8(1999), 1--14.

\item {[15]} G. Huisken and  C. Sinestrari,
             Convexity estimates for mean curvature flow and
             singularities of mean convex surfaces.
             Acta Math. 183(1999), 45-70.

\item {[16]}  N. Ivochkina, N.S. Trudinger, and X.-J. Wang,
              The Dirichlet problem for degenerate Hessian equations,
              Comm. Partial Diff. Eqns,  29(2004), 219-235.

\item {[17]} B. Kawohl,
       Rearrangements and convexity of level sets in PDE,
       Lecture Notes in Math., 1150.

\item {[18]}    N.V. Krylov,
       Nonlinear elliptic and parabolic equations of
       second order, Reidel, 1987.

\item {[19]}    N.V. Krylov,
       Smoothness of the payoff function for a controllable process
       in a domain, Math. USSR-Izv., 34(1990), 65-95.

\item {[20]}  N.V. Krylov,
        Lectures on fully nonlinear elliptic equations,
        Lipschitz Lectures, Univ. of Bonn, 1994.

\item {[21]}  D.G. Miguel,
       Differentiation of integrals in $\R^n$,
       Lecture Notes Math., 481 (1977).

\item {[22]} G. Perelman,
       The entropy formula for the Ricci flow and its geometric
       applications, \newline arXiv:math.DG/0211159.

\item {[23]} G. Perelman,
       Ricci flow with surgery on three manifolds,
       arXiv:math.DG/0303109.

\item {[24]} L. Simon,
        The minimal surface equation, Geometry V,
        Encyclopaedia Math. Sci., 90(1997), 239--272.

\item {[25]} N.S. Trudinger,
       On the Dirichlet problem for Hessian equations,
       Acta Math., 175(1995), 151-164.

\item {[26]} N.S. Trudinger and X.-J. Wang,
       The Bernstein problem for affine maximal hypersurfaces,
       Invent. Math., 140(2000), 399--422.

\item {[27]} B. White,
        The size of the singular set in mean curvature flow of
        mean-convex sets.
        J. Amer. Math. Soc. 13(2000), 665--695.

\item {[28]} B. White,
       The nature of singularities in mean curvature flow
       of mean-convex sets,
       J. Amer. Math. Soc., 16(2003), 123-138.

\item {[29]} G. Huisken and  C. Sinestrari,
       Mean curvature flow with serguries,
       Invent. Math., to appear.

\item {[30]} W. Sheng and X.J. Wang,
       Singularity profile in the mean curvature flow,
       Methods and Applications of Analysis, to appear.

\endRefs

\enddocument

\end

\enddocument

\end

\enddocument